\documentclass[aos,MSNbibl,seceqn,nameyear,dvips]{arximspdf}
\usepackage{mathbh}
\usepackage{accents}
\usepackage{graphicx}

\doi{10.1214/13-AOS1132} 
\volume{41}
\issue{4}
\pubyear{2013}
\firstpage{1970}
\lastpage{1998}

\makeatletter

\renewcommand{\mathring}[1]{\accentset{\circ}{#1}}

\newcommand{\rrVert}{\Vert}
\newcommand{\rrvert}{\vert}
\newcommand{\llVert}{\Vert}
\newcommand{\llvert}{\vert}

\newcommand{\cal}{\mathcal}

\newcommand{\be}{\mathbf{e}}
\newcommand{\btau}{\bolds{\tau}}
\newcommand{\bX}{\mathbf{X}}
\newcommand{\bzero}{\mathbf{0}}

\def\argmin{\mathop{\arg\min}}

\newtheorem{theo}{Theorem}[section]
\newtheorem{lem}{Lemma}[section]
\newtheorem{prop}{Proposition}[section]
\newtheorem{cor}{Corollary}[section]

\makeatother

\begin{document}
\begin{frontmatter}

\title{Inference in nonstationary asymmetric GARCH~models}
\runtitle{Inference in nonstationary asymmetric GARCH}

\begin{aug}
\author[A]{\fnms{Christian} \snm{Francq}\ead[label=e1]{Christian.Francq@ensae.fr}}
\and
\author[A]{\fnms{Jean-Michel} \snm{Zako\"{i}an}\corref{}\thanksref{t2}\ead[label=e2]{zakoian@ensae.fr}}
\runauthor{C. Francq and J.-M. Zako\"{i}an}
\affiliation{CREST and University Lille 3 (EQUIPPE)}
\address[A]{CREST\\
15 Boulevard Gabriel P\'{e}ri\\
92245 Malakoff cedex\\
France\\
\printead{e1}\\
\hphantom{E-mail: }\printead*{e2}}
\end{aug}

\thankstext{t2}{Supported by the Project ECONOM\&RISK (ANR 2010
blanc 1804 03).}

\received{\smonth{3} \syear{2013}}
\revised{\smonth{5} \syear{2013}}

%
\begin{abstract}
This paper considers the statistical inference of the class of
asymmetric power-transformed $\operatorname{GARCH}(1,1)$ models in presence of
possible explosiveness. We study the explosive behavior of
volatility when the strict stationarity condition is not met. This
allows us to establish the asymptotic normality of the quasi-maximum
likelihood estimator (QMLE) of the parameter, including the power
but without the intercept, when strict stationarity does not hold.
Two important issues can be tested in this framework: asymmetry and
stationarity. The tests exploit the existence of a universal
estimator of the asymptotic covariance matrix of the QMLE. By
establishing the local asymptotic normality (LAN) property in this
nonstationary framework, we can also study optimality issues.
\end{abstract}

%

\end{frontmatter}

\section{Introduction}

Following more than twenty years of tremendous development of the
theory of unit roots in linear time series models [see the seminal
papers by \citet{DicFul79} and \citet{PhiPer88}], there has
been, in the last decade, much interest in the statistical analysis of
nonlinear time series models under nonstationarity assumptions; see,
for example, \citet{KarTjs01}, \citet{KarMykTjs07},
\citet{LinLi08}, \citet{AueHor11}. In the framework of
GARCH (Generalized Autoregressive Conditional Heteroscedasticity)
models,
Jensen and Rahbek (\citeyear{JenRah04N1}, \citeyear{JenRah04N2}) were the first to establish an
asymptotic theory for the quasi-maximum likelihood estimator (QMLE)
of nonstationary $\operatorname{GARCH}(1,1)$, assuming that the intercept is fixed
to an arbitrary value. 
\citet{AknAlEHme11}
and \citet{AknAlE12} studied the properties of
weighted least-squares estimators.
\citet{FraZak12}
established the
asymptotic properties of the standard QMLE of the complete parameter
vector: they showed that, while the intercept cannot be consistently
estimated, the QMLE of the remaining parameters is consistent (in
the weak sense at the frontier of the stationarity region, and in
the strong sense outside) and asymptotically normal with or without
strict stationarity.
Asymptotic results for stationary $\operatorname{GARCH}(p,q)$ had been
established for the first time under mild conditions by \citet{BerHorKok03}.

Financial series are well known to present conditional asymmetry
features, in the sense that large negative returns tend to have more
impact on future volatilities than large positive returns of the same
magnitude. This stylized fact, known as the leverage effect, was first
documented by Black (\citeyear{Bla76}) and led to various generalizations of the
GARCH models of the first generation; see among others, \citet{GloJagRun93},
\citet{RabZak93},
\citet{HigBer92}, \citet{LiLi96}, \citet{FraZak10}.
Motivated by the Box--Cox
transformation, 
\citet{HwaKim04} introduced a power transformed ARCH model, and
the GARCH extension was studied by \citet{PanWanTon08}. In
this paper we consider an asymmetric power-transformed $\operatorname{GARCH}(1,1)$
model defined, for a given positive constant $\delta$, by
%
\begin{equation}
\label{ARCH1deJensenRahbek} \cases{\epsilon_t=h_t^{1/\delta}
\eta_t,
\cr
h_t=\omega_0+
\alpha_{0+}\bigl(\epsilon_{t-1}^+\bigr)^{\delta} +
\alpha_{0-}\bigl(-\epsilon_{t-1}^-\bigr)^{\delta} +
\beta_{0}h_{t-1},}
\end{equation}
with initial values $\epsilon_0$ and $h_0\geq0$, where 
$\omega_0>0$, $\alpha_{0+}\geq0$, $\alpha_{0-}\geq0$,
$\beta_{0}\geq0$, and using the notation $x^+=\max(x,0),
x^-=\min(x,0)$. In this model, $(\eta_t)$ is a sequence of
independent and identically distributed (i.i.d.) variables such that
%
\begin{equation}
\label{eta} 
E\eta_1^2=1 \quad\mbox{and}\quad P\bigl(
\eta_1^2=1\bigr)<1.
\end{equation}
Most commonly used extensions of the standard GARCH
of \citet{Eng82} and \citet{Bol86} can be written in the form
(\ref{ARCH1deJensenRahbek}).

The first goal of the present paper is to derive a strict
stationarity test in the framework of model
(\ref{ARCH1deJensenRahbek}). In this model, strict stationarity is
characterized by the negativity of the so-called top Lyapunov
exponent [see \citet{BouPic92}] which depends on the
parameters (except $\omega$) and the errors distribution. By
deriving
the asymptotic behavior of the QMLE of the
top-Lyapunov exponent, under stationarity and nonstationarity, a
strict stationarity test can be derived. 
The second goal of the paper is to propose a test for the symmetry
assumption in model~(\ref{ARCH1deJensenRahbek}), namely
$\alpha_{0+}=\alpha_{0-}$. Existing tests, to our knowledge, rely on
the stationarity assumption. Our aim is to derive a test which can be
used without bothering about stationarity.

The rest of the paper is organized as follows. In Section~\ref{Sec2}, we study the convergence of the volatility to infinity,
in a model encompassing (\ref{ARCH1deJensenRahbek}), when
stationarity does not hold. Section~\ref{QMLE} is devoted to the
asymptotic properties of the QMLE.
In Section~\ref{Sectest}, we consider strict stationarity testing and
asymmetry testing.
In Section~\ref{secLAN}, the LAN property is established and used to
derive the local asymptotic power of the proposed tests. Local
alternative allowing for an arbitrary rate of convergence with
respect to $\omega_0$ are considered.
Optimality
issues are discussed. Necessary and sufficient conditions on the
noise density are derived for the tests to be uniformly locally
asymptotically most powerful. Section~\ref{Sec4} is devoted to the
case where the power $\delta$ is unknown and is jointly estimated
with the volatility coefficients.
Proofs and 
technical lemmas are in Section~\ref{complement}. The possibility of
extensions is discussed in Section~\ref{conclu}. Due to space restrictions, several lemmas and proofs,
along with a
study of the finite sample performance of the
stationarity and asymmetry tests and an empirical application, are
included in the supplementary file [\citet{FraZak}].

\section{Explosivity in the augmented $\operatorname{GARCH}(1,1)$}\label{Sec2}

In this section, we analyze the convergence of the volatility to
infinity, for a class of augmented GARCH processes encompassing
(\ref{ARCH1deJensenRahbek}) and many $\operatorname{GARCH}(1,1)$ models
introduced in the literature; see \citet{Hor08}. Given a sequence
$(\xi _t)_{t\geq 0}$, let $(\epsilon_t)_{t\geq1}$ be defined by
%
\begin{equation}
\label{modgarch11augm} \cases{ \epsilon_t= h_t^{1/\delta}
\xi_t, \qquad t=1,2,\ldots,
\cr
h_t=\omega(
\xi_{t-1})+ a(\xi_{t-1})h_{t-1},}
\end{equation}
where $\delta$ is a positive constant, $h_0\geq0$ is a given
initial value and the functions $\omega(\cdot)$ and $a(\cdot)$
satisfy $\omega\dvtx  \mathbb{R} \to[\underline{\omega}, +\infty)$ and
$a\dvtx  \mathbb{R} \to[0,+\infty)$, for some $\underline{\omega}>0$.
When $(\xi_t)$ is assumed to be a white noise, $(\epsilon_t)$ is
called an augmented GARCH process. We purposely use a different
notation for $\xi_t$ in (\ref{modgarch11augm}) and $\eta_t$
in (\ref{ARCH1deJensenRahbek}) because, for the moment, we only assume
that $(\xi_t)$ is stationary and ergodic. Define in
$\mathbb{R}\cup\{+\infty\}$ the top Lyapunov exponent
\[
\gamma=E\log a(\xi_1).
\]
The following proposition is an extension of results proven for the
standard $\operatorname{GARCH}(1,1)$ by \citet{Nel90} and completed by
\citet{KluLinMal04} and \citet{FraZak12}.
%
\begin{prop}\label{convtoinfty}
For the process $(\epsilon_t)$ satisfying (\ref{modgarch11augm}),
the following properties hold:
\begin{enumerate}[(ii)]
\item[(i)] When $\gamma>0$, $ h_t\to\infty$ a.s. at an
exponential rate: for any $\rho> e^{-\gamma}$,
\[
\rho^th_t \to\infty \quad\mbox{and}\quad \mbox{if $E\bigl|\log\bigl(
\xi_1^2\bigr)\bigr|<\infty$}\qquad \rho ^t
\epsilon_t^2 \to\infty\qquad\mbox{a.s.} \mbox{ as } t\to\infty.
\]
\item[(ii)]
When $\gamma=0$ and $(\xi_t)$ is time reversible [i.e., for
all $k$ the distributions of $(\xi_t,\xi_{t-1},\ldots, \xi_{t-k})$ and
$(\xi_{t-k},\ldots,\xi_{t-1}, \xi_t)$ are identical], the following
convergences in probability hold as $t\to\infty$:
\[
h_t \to\infty\quad\mbox{and}\quad \mbox{if } E\bigl|\log\bigl(\xi_1^2
\bigr)\bigr|<\infty\qquad \epsilon^2_t \to\infty.
\]
Moreover, if $\psi$ is a decreasing bijection from $(0,\infty)$
to $(0,\infty)$, if\break $E\psi(h_1)<\infty$ [resp.,
$E\psi(\epsilon_1^2)<\infty$ and $E|{\log(\xi_1^2)}|<\infty$], then
%
\begin{equation}
\label{empmeans} \psi(h_t)\to0 \qquad\bigl[\mbox{resp., }\psi\bigl(
\epsilon_t^2\bigr)\to0\bigr]\qquad\mbox{in
$L^1$.}
\end{equation}
\end{enumerate}
\end{prop}
The main ideas of the proof are as follows. The a.s. convergence of
$h_t$ to infinity in the case $\gamma
>0$ follows from the minoration
$\log h_t\geq\log\underline{\omega}+\sum_{i=1}^{t-1}\log
a(\xi_{t-i})$, and the fact that the latter sum is strictly
increasing, in average, as $t$ goes to infinity. The argument is in
failure when $\gamma=0$, the expectation of the sum being equal to
zero.
The key argument in this case is that the sequence $(h_t)$ is
increasing \textit{in distribution}. Indeed, taking $h_0=0$ we have
$h_1=\omega(\xi_0)$ and $h_2=\omega(\xi_1)+
a(\xi_0)\omega(\xi_0)\stackrel{d}{=}\omega(\xi_0)+
a(\xi_1)\omega(\xi_1)>h_1$ under the reversibility assumption, and
the same argument applies for any
$t>0$.

In the rest of the paper, these results will be applied with
$\xi_t=\eta_t$ to model (\ref{ARCH1deJensenRahbek}), for which the
top Lyapunov exponent
is given by
%
\[
\gamma_0=E\log a_0(\eta_1),\qquad
a_0(x)=\alpha_{0+}\bigl(x^+\bigr)^{\delta
}+
\alpha_{0-}\bigl(-x^-\bigr)^{\delta}+\beta_0.
\]

\section{Asymptotic properties of the QMLE}
\label{QMLE}
We wish to estimate $\vartheta_0=(\alpha_{0+},\break
\alpha_{0-},\beta_0)'$ from observations $\epsilon_t, t=1,\ldots,
n$, in the stationary and the explosive cases under mild assumption.
Denote by $\theta=(\omega,\alpha_+,\alpha_-,
\beta)'$ the parameter and define the QMLE as any
measurable solution of
%
\begin{eqnarray}
\label{estimateurQMVJensenRahbek} \hat{\theta}_n&=&(\hat{\omega
}_n,\hat{\alpha}_{n+}, \hat{\alpha}_{n-}, \hat{
\beta}_n)'= \argmin_{\theta\in
\Theta}
\frac{1}{n}\sum_{t=1}^n
\ell_t(\theta),\nonumber\\[-8pt]\\[-8pt]
\ell_t(\theta)&=&\frac{\epsilon_t^2}{\sigma_t^2(\theta)}+\log
\sigma_t^2(\theta),\nonumber
\end{eqnarray}
where $\Theta$ is a compact subset of $(0,\infty)^4$ containing the
true value $\theta_0=(\omega_0,\alpha_{0+}, \alpha_{0-},\beta
_0)'$, and 
$\sigma_t^{\delta}(\theta)=\omega+\alpha_+
(\epsilon_{t-1}^+)^{\delta}+\alpha_-
(-\epsilon_{t-1}^-)^{\delta}+\beta\sigma_{t-1}^{\delta}(\theta)$ for
$t=1,\ldots,n$ [with initial values for $\epsilon_0$ and
$\sigma_0^{\delta}(\theta)$]. The rescaled residuals are defined by
$\hat{\eta}_t=\eta_t(\hat{\theta}_n)$ where
$\eta_t(\theta)=\epsilon_t/\sigma_t(\theta)$ for $t=1,\ldots,n$.

Write $\vartheta=(\alpha_{+}, \alpha_{-}, \beta)'$ and let
$\hat{\vartheta}_n= (\hat{\alpha}_{n+}, \hat{\alpha}_{n-},
\hat{\beta}_n )'$.

\subsection{\texorpdfstring{Consistency and asymptotic normality of $\hat{\vartheta}_n$}
{Consistency and asymptotic normality of theta n}}

The following theorem extends, to the nonstationary framework, results
obtained for the stationary
case [see \citet{HamZak11} and the references therein],
which we recall for convenience. 
We introduce
the assumptions:

{A1:}
The support of $(\eta_t)$ contains at least 3 points and is not
concentrated on the positive or the negative line.

{A2:} When $t$ tends to infinity,
\[
E \Biggl\{1+\sum_{i=1}^{t-1}a_0(
\eta_{1})\cdots a_0(\eta_{i}) \Biggr\}
^{-1}= o \biggl(\frac{1}{\sqrt{t}} \biggr).
\]
Note that A2, which is only required in the case $\gamma_0=0$,
is obviously satisfied in the degenerate case when
$a(\eta_t)=1$, a.s.,
since the expectation is then equal to $1/t$.

To handle
initial values we introduce the following notation. For any
asymptotically stationary process $(X_t)_{t\geq0}$, let
$E_{\infty}(X_t)=\lim_{t\to\infty} E(X_t)$ provided this limit
exists. 
Let also $\mathring{\Theta}$ denote the interior of $\Theta$.
%
\begin{theo}
\label{theoconsistenceJensenRahbek}
Let (\ref{ARCH1deJensenRahbek})--(\ref{eta}) and \textup{A1} hold. Then
the QMLE defined in (\ref{estimateurQMVJensenRahbek}) satisfies the
following properties:
%
\begin{enumerate}[(iii)]
\item[(i)] \textup{Stationary case.} When $\gamma_0< 0$, and $\beta<1$
for all
$\theta\in\Theta$,
\[
\hat{\theta}_n\to\theta_0\qquad \mbox{a.s. as } n\to\infty.
\]
If, in addition, $\kappa_{\eta}=E\eta_1^4 \in(1,\infty)$ and
$\theta_0 \in\mathring{\Theta}$,
we have
%
\begin{equation}
\label{normalityQMVusuel} \sqrt{n} (\hat{\theta}_n-
\theta_0 )\stackrel{d} {\to} {\cal N} \bigl\{0,(\kappa_{\eta}-1){
\cal J}^{-1} \bigr\}\qquad\mbox{as } n\to\infty,
\end{equation}
where
%
\begin{equation}
\label{jpourungarch11} {\cal J} =\frac{4}{\delta^{2}} E_{\infty} \biggl(
\frac{1}{\sigma_t^{2\delta}}\,\frac{\partial\sigma
_t^\delta}{\partial
\theta}\,\frac{\partial\sigma_t^\delta}{\partial
\theta'}(\theta_0)
\biggr).
\end{equation}
\item[(ii)] \textup{Explosive case.}
When $\gamma_0>0$, if $P(\eta_1=0)=0$,
\[
\hat{\vartheta}_n\to \vartheta_0\qquad \mbox{a.s. as } n\to
\infty.
\]
If, in addition, $\kappa_{\eta}\in(1,\infty)$,
$E|{\log\eta_1^2}|<\infty$ and $\theta_0\in
\mathring{\Theta}$,
%
\begin{equation}
\label{normalitydeleQMVdealpha} \sqrt{n} (\hat{\vartheta}_n-
\vartheta_0 )\stackrel {d} {\to} {\cal N} \bigl\{0,(
\kappa_{\eta}-1){\cal I}^{-1} \bigr\} 
\end{equation}
as $n\to\infty$, where ${\cal I}$ is a positive definite matrix.
%
\item[(iii)] \textup{At the boundary of the stationarity region.}
When $\gamma_0= 0$, if $P(\eta_1=0)=0$, and
$\forall\theta\in\Theta$, $\beta<\|1/a_0(\eta_1)\|_{p}^{-1}$
for some $p>1$,
%
\[
\hat{\vartheta}_n\to \vartheta_0\qquad \mbox{in probability
as } n\to\infty.
\]
If, in addition,
$\theta_0\in \mathring{\Theta}$,
$\kappa_{\eta}\in(1,\infty)$, $E|{\log\eta_1^2}|<\infty$
and \textup{A2} is satisfied,
then (\ref{normalitydeleQMVdealpha}) holds.
\end{enumerate}
\end{theo}
%
The key ideas of the proof can be summarized as follows. First, we
note that $\hat{\theta}_n$ can be equivalently defined as the
minimizer of $\frac{1}{n}\sum_{t=1}^n
\{\ell_t(\theta)-\ell_t(\theta_0)\}$, where
$\ell_t(\theta)-\ell_t(\theta_0)$ is a function of $\eta_t^2$ and
the ratio $\sigma_t^{\delta}(\theta)/h_t$. While the numerator and
the denominator explode to infinity as $t$ increases, the ratio is
close to a stationary process for $t$ sufficiently large. For
instance, in the symmetric ARCH(1) case ($\alpha_+=\alpha_-=\alpha$
and $\beta=0$), we have $\sigma_t^{\delta}(\theta)/h_t\to
\alpha/\alpha_0$, a.s. in the strictly explosive case (in
probability in the case $\gamma=0$). The situation is much more
intricate when $\beta\ne0$, but we can show that, when $\gamma>0$,
\[
\biggl\llvert \frac{\sigma_t^{\delta}(\theta)}{h_t}-v_{t}(\vartheta )\biggr\rrvert \to0
\qquad\mbox{a.s. as $t\to\infty$}
\]
uniformly on some compact set included in $\Theta$, where
$(v_{t}(\vartheta))$ is a strictly stationary and ergodic process. The
a.s. convergence is replaced by a $L^p$ convergence in the case
$\gamma=0$. The consistency results are established by showing that the
criterion in which $\sigma_t^{\delta}(\theta)/h_t$ is replaced by
$v_{t}(\vartheta)$ produces an estimator which is consistent
to~$\vartheta_0$. Similar arguments are used to prove the asymptotic
normality results, but we now show that
\[
\biggl\llVert \frac{1}{\sigma^{\delta}_t(\theta)}\,\frac{\partial
\sigma_t^{\delta}}{\partial\vartheta}(\theta_0)-
d_t\biggr\rrVert \to 0\qquad\mbox{in $L^p$ as $t\to\infty$}
\]
for some strictly stationary and ergodic process $d_t$.

An explicit expression of ${\cal I}$ is given in the supplementary
file [\citet{FraZak}]. To conclude the section, it can be
noted that no asymptotically valid inference on $\omega_0$ can be
done in the nonstationary case; see Propositions 2.1 and 3.1 in
\citet{FraZak12}, denoted hereafter FZ, for the standard
$\operatorname{GARCH}(1,1)$ model.

\subsection{\texorpdfstring{A universal estimator of the asymptotic variance of $\hat{\vartheta}_n$}
{A universal estimator of the asymptotic variance of theta n}}

In view of (\ref{normalityQMVusuel})--(\ref{jpourungarch11}), when
$\gamma_0<0$ the asymptotic distribution of the QMLE
$\hat{\vartheta}_n$ of $\vartheta_0$ (the parameter without $\omega
_0$) is given by
%
\begin{equation}
\label{normalityQMVusueldealpha} \sqrt{n} (\hat{\vartheta}_n-
\vartheta_0 )\stackrel {d} {\to} {\cal N} \bigl\{0,(
\kappa_{\eta}-1){\cal I}_*^{-1} \bigr\} \qquad\mbox{as } n\to
\infty
\end{equation}
with
%
\begin{equation}
\label{upsilon} {\cal I}_* ={\cal J}_{\vartheta, \vartheta}-{\cal J}_{\vartheta,
\omega}{\cal
J}_{\omega, \omega}^{-1}{\cal J}_{\omega, \vartheta}, 
\end{equation}
${\cal J}_{\omega, \omega}=
\frac{4}{\delta^{2}}E_{\infty} (\frac1{h_t^2}\,\frac{\partial
\sigma_t^{\delta}}{\partial\omega}\,\frac{\partial
\sigma_t^{\delta}}{\partial\omega}(\theta_0) ), {\cal
J}_{\vartheta, \vartheta}=
\frac{4}{\delta^{2}}E_{\infty} (\frac1{h_t^2}\,\frac{\partial
\sigma_t^{\delta}}{\partial\vartheta}\,\frac{\partial
\sigma_t^{\delta}}{\partial\vartheta'}(\theta_0) )$ and
${\cal
J}_{\omega, \vartheta}= {\cal J}_{\vartheta, \omega}'=
\frac{4}{\delta^{2}}E_{\infty} (\frac1{h_t^2}\,\frac{\partial
\sigma_t^{\delta}}{\partial\omega}\,\frac{\partial
\sigma_t^{\delta}}{\partial\vartheta'}(\theta_0) )$.
Letting
\[
\hat{{\cal J}}_{\vartheta,
\vartheta}=\frac{4}{\delta^{2}}\frac1n \sum
_{t=1}^n \frac1{\sigma_t^{2\delta}(
\hat{\theta}_n)}\,\frac{\partial
\sigma_t^{\delta}}{\partial\vartheta}\,\frac{\partial
\sigma_t^{\delta}}{\partial\vartheta'}(\hat{
\theta}_n)
\]
and
defining $\hat{{\cal J}}_{\vartheta, \omega}, \hat{{\cal
J}}_{\omega, \omega}$ and $\hat{{\cal J}}_{\omega, \vartheta}$
accordingly, it can be shown that
\[
\hat{\cal I}_* =\hat{{\cal J}}_{\vartheta, \vartheta
}-
\hat{{\cal J}}_{\vartheta, \omega}\hat{{\cal J}}_{\omega, \omega}^{-1}
\hat{{\cal J}}_{\omega, \vartheta}
\]
is a strongly consistent estimator of ${\cal I}_*$ in the stationary
case $\gamma_0<0$.
The following result shows that this estimator also provides
a consistent estimator of the asymptotic variance of
$\hat{\vartheta}_n$
in the nonstationary case $\gamma_0 \geq0$.
%
\begin{theo}
\label{propo1} Let the assumptions required for the consistency
results in Theorem
\ref{theoconsistenceJensenRahbek} hold, assume ${\kappa}_{\eta}\in
(1, \infty)$ and let
$\hat{\kappa}_{\eta}=n^{-1}\sum_{t=1}^n\hat{\eta}_t^4$, where
$\hat{\eta}_t=\epsilon_t/\sigma_t(\hat{\theta}_n)$.
\begin{enumerate}[(iii)]
\item[(i)] When $\gamma_0<0$, we have $\hat{\kappa}_{\eta}\to
{\kappa}_{\eta}$ and $\hat{\cal I}_* \to{\cal I}_*$ a.s. as
$n\to\infty$.
\item[(ii)] When $\gamma_0>0$, we have $\hat{\kappa}_{\eta}\to
{\kappa}_{\eta}$ and $\hat{\cal I}_*\to{\cal I}$ a.s.
\item[(iii)] When $\gamma_0=0$, we have $\hat{\kappa}_{\eta}\to
{\kappa}_{\eta}$ and, if \textup{A2} is satisfied, $\hat{\cal
I}_*\to{\cal I}$ in probability.
\end{enumerate}
In any case, $(\hat{\kappa}_{\eta}-1)\hat{\cal I}_*^{-1}$ is a
consistent estimator of the asymptotic variance of the QMLE of
$\vartheta_0$.
\end{theo}

It follows that asymptotically valid confidence intervals for the
parameter $\vartheta_0$
can be constructed without knowing if the underlying process is
stationary or not. This theorem also has interesting applications for
testing problems,
which we now consider.

\section{Testing}
\label{Sectest} In this section we consider testing stationarity and
testing asymmetry.
\subsection{Strict stationarity testing}
\label{S5}

Consider the strict stationarity testing problems
%
\begin{equation}
\label{testing^problem} H_0\dvtx  \gamma_0<0
\quad\mbox{against}\quad H_1\dvtx  \gamma_0\geq0
\end{equation}
and
%
\begin{equation}
\label{testing^problem2} H_0\dvtx  \gamma_0
\geq0 \quad\mbox{against}\quad H_1\dvtx  \gamma_0< 0.
\end{equation}
Let
$\hat{\gamma}_n={\gamma}_n(\hat{\theta}_n)$ be the empirical
estimator of $\gamma_0$, with for any $\theta\in\Theta$,
%
\begin{equation}
\label{gammahat} \gamma_n(\theta)=\frac1n\sum
_{t=1}^n\log \bigl[\alpha_{+} \bigl\{
\eta^+_t(\theta) \bigr\}^{\delta}+\alpha_{-} \bigl\{-
\eta ^-_t(\theta) \bigr\}^{\delta}+\beta \bigr],
\end{equation}
where $\eta_t(\theta)=\epsilon_t/\sigma_t(\theta)$. The following
result shows that the asymptotic distribution of $\hat{\gamma}_n$
is particularly simple in the nonstationarity case.
%
\begin{theo}
\label{jointdistribution}
Let $u_t=\log a_0(\eta_t)-\gamma_0$, and $\sigma^2_u=Eu_t^2$. Then,
under the assumptions of Theorem \ref{theoconsistenceJensenRahbek},
%
\begin{equation}
\label{amontrerenappendix3} \sqrt{n}(\hat{\gamma}_n-
\gamma_0)\stackrel{d} {\to} {\cal N} \bigl(0,\sigma_{\gamma}^2
\bigr) \qquad\mbox{as }n\to\infty,\vadjust{\goodbreak}
\end{equation}
where
\[
\sigma_{\gamma}^2=\cases{ \sigma_u^2+
(\kappa_{\eta}-1) 
\bigl\{a' {\cal
J}^{-1}a-(1-\nu_1)^2\bigr\}, &\quad when $
\gamma_0<0$,
\vspace*{2pt}\cr
\sigma_u^2, &\quad when $
\gamma_0\geq0$,}
\]
with $a=(0, \tilde{\nu}_{1,+}, \tilde{\nu}_{1,-}, \nu_1/\beta
_0)'$ and
\[
\tilde{\nu}_{1+}=E \biggl\{\frac{(\eta^+_1)^{\delta}}{a_0(\eta
_1)} \biggr\},\qquad \tilde{
\nu}_{1-}=E \biggl\{\frac{(-\eta
^-_1)^{\delta}}{a_0(\eta_1)} \biggr\},\qquad \nu_1=E
\biggl\{\frac
{\beta_0}{a_0(\eta_1)} \biggr\}.
\]
\end{theo}

Let $\hat{\sigma}_u^2$ be the empirical variance of $\log
\{\hat{\alpha}_{n+} (\hat{\eta}^+_t )^{\delta
}+\hat{\alpha}_{n-} (-\hat{\eta}^-_t )^{\delta}+\hat
{\beta}_n \}$,
for $t=1,\ldots,n$. Under the assumptions of
Theorem \ref{jointdistribution}, it can be shown that
$\hat{\sigma}_u^2$ is a weakly consistent estimator of $\sigma_u^2$.
The statistics
\[
T_n=\sqrt{n} {\hat{\gamma}_n}/{\hat{
\sigma}_u}
\]
are
thus asymptotically ${\cal N} (0,1)$ distributed when $\gamma_0=0$.
For the testing
problem (\ref{testing^problem}) [resp., (\ref{testing^problem2})], at
the asymptotic significance level $\underline{\alpha}$, this leads
to consider the critical region
%
\begin{equation}
\label{test} {\mathrm C}^{\mathrm{ST}}= \bigl\{T_n>
\Phi^{-1}(1-\underline{\alpha }) \bigr\} \qquad\bigl[\mbox{resp., ${\mathrm
C}^{\mathrm{NS}}= \bigl\{ T_n<\Phi^{-1}(\underline{
\alpha}) \bigr\}$}\bigr].
\end{equation}

\subsection{Asymmetry testing}
\label{Sasymmetry} It is of particular interest to test the
existence of a leverage effect in stock market returns. In the
framework of model (\ref{ARCH1deJensenRahbek}), this testing problem
is of the form
%
\begin{equation}
\label{testingasymmetry} H_0\dvtx  \alpha_{0+}=
\alpha_{0-} \quad\mbox{against}\quad H_1\dvtx  \alpha_{0+}\neq
\alpha_{0-}.
\end{equation}
Consider the test statistic for symmetry
\[
T_n^{\mathrm{S}}:=\frac{\sqrt{n}(\hat{\alpha}_{n+}-\hat{\alpha
}_{n-})}{\hat{\sigma}_{T^\mathrm{S}}},\qquad \hat{\sigma}_{T^\mathrm{S}}=
\sqrt{(\hat{\kappa}_{\eta}-1)\be '\hat{\cal
I}_*^{-1}\be}
\]
with $\be'=(1,-1,0)$. The following result is a direct consequence
of (\ref{normalitydeleQMVdealpha}), (\ref{normalityQMVusueldealpha})
and Theorem \ref{theoconsistenceJensenRahbek}.
%
\begin{cor}
\label{coro31} Assume that $\theta_0 \in
\mathring{\Theta}$ and the assumptions of Theorem~\ref{theoconsistenceJensenRahbek}
 hold. For the testing problem (\ref{testingasymmetry}),
the test defined by the critical region
%
\begin{equation}
\label{test3} {\mathrm C}^{\mathrm{S}}= \bigl\{\bigl|T_n^{\mathrm{S}}\bigr|>
\Phi ^{-1}(1-\underline{\alpha}/2) \bigr\}
\end{equation}
has the asymptotic significance level $\underline{\alpha}$ and is
consistent.
\end{cor}
We emphasize the fact that this test for symmetry does not require
any stationarity assumption. The somewhat surprising output is that
the usual Wald test, based on the asymptotic theory for the
stationary case, also works in the nonstationary
situation.\setcounter{footnote}{1}\footnote{For instance, in ARMA models, Wald tests on the
parameters are not the same in the stationary and nonstationary
cases.}

\section{Asymptotic local powers}
\label{secLAN} This section investigates the asymptotic behavior under
local alternatives of the asymmetry test (\ref{test3}) and of the
strict stationarity test (\ref{test}). We first establish the LAN of
the power-transformed GARCH model without imposing any stationarity
constraint. This LAN property will be used to derive the asymptotic
properties of our tests, but the result is of independent interest; see
\citet{van98} for a general reference on LAN and its applications,
and see \citet{DroKla97}, \citet{DroKlaWer97} and
\citet{LinMcA03} for applications to GARCH and other stationary
processes.

\subsection{LAN without stationarity constraint}
Assume that $\eta_t$ has a density $f$ which is positive everywhere,
with third-order derivatives such that
%
\begin{equation}\label{hypo2surf}
\lim_{|y|\to\infty}yf(y)=0 \quad\mbox{and}\quad \lim_{|y|\to
\infty}y^2f'(y)=0,
\end{equation}
and that, for some positive constants $K$ and $\delta$,
%
\begin{eqnarray}\label{hypo3surf}
|y|\biggl\llvert \frac{f'}{f}(y)\biggr\rrvert +y^2\biggl
\llvert \biggl(\frac{f'}{f} \biggr)'(y)\biggr\rrvert
+y^2\biggl\llvert \biggl(\frac{f'}{f} \biggr)''(y)
\biggr\rrvert &\leq& K \bigl(1+|y|^{\delta} \bigr),
\\
\label{hypo3surfbis}
E\llvert \eta_1\rrvert ^{2\delta}&<&\infty.
\end{eqnarray}
These regularity conditions are satisfied for numerous
distributions, in particular for the Gaussian distribution with
$\delta=2$, and entail the existence of the Fisher information for
scale
\[
\iota_{f}=\int \bigl\{1+yf'(y)/f(y) \bigr
\}^2f(y)\,dy<\infty.
\]
Given the initial values $\epsilon_0$ and $h_0$, the density of the
observations $(\epsilon_1,\ldots, \epsilon_n)$ satisfying
(\ref{ARCH1deJensenRahbek}) is given by
$L_{n,f}(\theta_0)=\prod_{t=1}^n \sigma_t^{-1}(\theta_0)
f \{\sigma_t^{-1}(\theta_0)\epsilon_t \}$. Around \mbox{$\theta_0
\in\mathring{\Theta}$}, let a sequence of local parameters of
the form
%
\begin{equation}
\label{altasusual} \theta_n=\theta_0+
\btau_n/\sqrt{n},
\end{equation}
where $(\btau_n)$ is a bounded
sequence of $\mathbb{R}^4$. Without loss of generality, assume that
$n$ is sufficiently large so that $\theta_n\in{\Theta}$. Under the
strict stationarity condition $\gamma_0<0$, \citet{DroKla97} showed that, for standard GARCH, the log-likelihood ratio
$\Lambda_{n,f}(\theta_n,\theta_0)=\log
L_{n,f}(\theta_n)/L_{n,f}(\theta_0)$ satisfies the LAN property
%
\begin{equation}
\label{lan} \Lambda_{n,f}(\theta_n,\theta_0)={
\btau_n'}S_{n,f}(\theta_0)-
\tfrac{1}{2}\btau_n'\mathfrak{I}_{f}
\btau_n+o_{P_{\theta_0}}(1),
\end{equation}
where $S_{n,f}(\theta_0)\stackrel{d}{\longrightarrow}{\cal
N} \{0,\mathfrak{I}_{f} \}$ under $P_{\theta_0}$ as
$n\to\infty$. Note that the so-called central sequence $S_{n,f}$ is
conditional on the initial values. In the stationary case, \citet{LeeTan05} showed that the initial values have no influence on
the LAN property. The following proposition shows that (\ref{lan})
holds regardless of $\gamma_0$.
%
\begin{prop}
\label{LANproperty} When $\theta_0 \in \mathring{\Theta}$,
under (\ref{hypo2surf})--(\ref{hypo3surfbis}) we have the LAN
property (\ref{lan}). When $\gamma_0< 0$, we have
$\mathfrak{J}_{f}=\frac{\iota_f}4{\cal J}$, where ${\cal J}$ is
defined in (\ref{jpourungarch11}). When $\gamma_0\geq0$, the Fisher
information is the degenerate matrix
%
\begin{equation}
\label{Infomat} 
\mathfrak{I}_{f}=\frac{\iota_{f}}{4} \pmatrix{0&0_3'
\cr
0_3&{\cal I}},
\end{equation}
where ${\cal I}$ is the positive definite matrix introduced in (\ref
{normalitydeleQMVdealpha}).
\end{prop}

\subsection{\texorpdfstring{Near-global alternatives with respect to $\omega_0$}
{Near-global alternatives with respect to omega 0}}
We now
show that, in the nonstationary case, LAN continues to hold when
the local alternative allows for an arbitrary rate of convergence
with respect to $\omega_0$. To this aim we assume that
%
\begin{equation}
\label{newalt} \theta_n=\theta_0+\upsilon_n
\be _1+\frac{\btau_n}{\sqrt{n}},
\end{equation}
where $\be_1=(1, 0, 0,
0)'$, $(\btau_n)$ is as in (\ref{altasusual}), and $(\upsilon_n)$ is
a deterministic sequence converging to zero.
The next result shows that, in
the nonstationary case, (\ref{lan}) which was established under
(\ref{altasusual}), continues to hold under the more general
alternatives (\ref{newalt}). For simplicity, take
$\btau_n=\btau=(\tau_1, \tilde{\btau}')'$ and
$\tilde{\btau}'=(\tau_2,\tau_3,\tau_4)$.
%
\begin{prop}
\label{LANpropertybis} Let $\theta_0 \in \mathring{\Theta}$
with $\gamma_0\geq0$. Then, under
(\ref{hypo2surf})--(\ref{hypo3surfbis}) and (\ref{newalt}), we have
the LAN property
\[
\Lambda_{n,f}(\theta_n,
\theta_0) \stackrel{d} {\longrightarrow} {\cal N} \biggl(-
\frac{\iota_{f}}{8}\tilde{\btau}'{\cal I}\tilde{\btau},
\frac{\iota_{f}}{4} \tilde{\btau}'{\cal I}\tilde{\btau} \biggr)
\qquad\mbox{under $P_{\theta_0}$ as $n\to\infty$}.
\]
\end{prop}
Note that this Gaussian law is the distribution of the
log-likelihood ratio in the statistical model ${\cal
N} \{\tilde{\btau},4{\cal I}^{-1}/\iota_{f} \}$ of
parameter $\tilde{\btau}$, or equivalently in the statistical model
${\cal N} \{\iota_{f}{\cal I}\tilde{\btau}/4,\iota_{f}{\cal
I}/4 \}$. To interpret this result in terms of convergence of
statistical experiments [see \citet{van98} for details],
assume that $\upsilon_n=\upsilon\nu_n$ where $\upsilon\in
\mathbb{R}$ and $(\nu_n)$ is a given sequence converging to zero as
$n\to\infty$. Denoting by ${\cal T}$ a subset of $\mathbb{R}^4$
containing a neighborhood of $\bzero$, the so-called local
experiments
$ \{L_{n,f}(\theta_0+\upsilon\nu_n\be_1+(0,\tilde{\btau
}')/\sqrt{n}),
(\upsilon, \tilde{\btau}')\in{\cal T} \}$ converge to the
Gaussian experiment $ \{{\cal N} (\tilde{\btau},4{\cal
I}^{-1}/\iota_{f} ),(\upsilon, \tilde{\btau}')\in{\cal
T} \}$.

Interestingly, the parameter $\upsilon$ vanishes in the limiting
experiment. Consequently, in the limit experiment there exists no
test on the parameter $\upsilon$ (except of trivial power equal to
the level). On the other hand, the limit of any converging sequence
of power functions in the local experiments is a power function in
the Gaussian limit experiment, by the asymptotic representation
theorem. We can conclude that there exists no test with a nontrivial
asymptotic power, for local alternatives on the parameter
$\upsilon$ at the rate $1/\nu_n$. Given that the rate of convergence
of $\nu_n$ to zero is arbitrary, the LAN approach shows that no
asymptotically valid inference can be made on the parameter
$\omega_0$.\footnote{This is in accordance with the observation
that, at least in the explosive case, the Fisher information with
respect to $\omega_0$ is bounded as $n$ increases. A proof is
available from the authors.}

\subsection{Local asymptotic powers of the tests}
The LAN property, with the help of Le Cam's third lemma, allows us
to easily compute local asymptotic powers of tests. In view
of Theorem \ref{jointdistribution},
\[
\lim_{n\to\infty}P_{\theta_0} \bigl({\mathrm C}^{\mathrm
{ST}}
\bigr)=\lim_{n\to\infty}P_{\theta_0} \bigl({\mathrm
C}^{\mathrm{NS}} \bigr)=\underline{\alpha},
\]
when $\theta_0$ is such that $\gamma_0=0$. For
$\btau$ 
such that
$\theta_0+\btau/\sqrt{n}\in\Theta$, we denote by $P_{n,\btau}$ the
distribution of the observations $(\epsilon_1,\ldots,\epsilon_n)$
when the parameter is $\theta_0+\btau/\sqrt{n}$. We should use the
notation $(\epsilon_{1,n},\ldots,\epsilon_{n,n})$ instead of
$(\epsilon_1,\ldots,\epsilon_n)$ because the parameter varies with
$n$, but we will avoid this heavy notation. Let
\[
a_{\btau}(\eta_1)= \biggl(\alpha_{0+}+
\frac{\tau_2}{\sqrt {n}} \biggr) \bigl(\eta_1^+ \bigr)^{\delta}+
\biggl(\alpha _{0-}+\frac{\tau_3}{\sqrt{n}} \biggr) \bigl(-
\eta_1^- \bigr)^{\delta}+\beta_0+
\frac{\tau_4}{\sqrt{n}}.
\]
Local alternatives for the ${\mathrm C}^{\mathrm{ST}}$-test (resp.,
the ${\mathrm C}^{\mathrm{NS}}$-test) are obtained for $\btau$ such
that $E\log a_{\btau}(\eta_1)>0$ (resp., $E\log
a_{\btau}(\eta_1)<0$).

\begin{prop}
\label{LAP2tests} Under the assumptions of Theorem
\ref{theoconsistenceJensenRahbek} and Proposition~\ref{LANproperty},
the local
asymptotic powers of the strict stationarity tests (\ref{test}) are
given by
%
\begin{equation}
\label{LPASTtest} \lim_{n \to
\infty}P_{n,\btau} \bigl({\mathrm
C}^{\mathrm{ST}} \bigr)=\Phi \bigl\{c_f(\theta_0)-
\Phi^{-1}(1-\underline{\alpha}) \bigr\}
\end{equation}
and, using the notation of Theorem
\ref{jointdistribution},
\[
\lim_{n \to
\infty}P_{n,\btau} \bigl({\mathrm C}^{\mathrm{NS}}
\bigr)=\Phi \bigl\{\Phi^{-1}(\underline{\alpha})-c_f(
\theta_0) \bigr\},
\]
where
\[
c_f(\theta_0)=\frac{ (\tau_2\tilde{\nu}_{1+}+\tau_3\tilde
{\nu}_{1-}+\tau_4\nu_1/\beta_0 )E\log a_0(\eta_1) \{
1+\eta_1{f'(\eta_1)}/{f(\eta_1)} \}}{\delta\sigma
_u(1-\nu_1)}.
\]
\end{prop}
We now compute the local asymptotic power of the asymmetry test
defined by (\ref{test3}). We thus consider a sequence of local
parameters of the form $\theta_n=\theta_0+\btau/\sqrt{n}$ where
$\theta_0=(\omega_0,\alpha_0,\alpha_0,\beta_0)'$ and
$\btau=(\tau_1,\tau_2,\tau_3,\tau_4)'$ (with $\tau_2\neq\tau_3$
under a local alternative). We denote by $P_{n,\btau}^{\mathrm{S}}$
the distribution of the observations under the assumption that
the parameter is $\theta_n$.
%
\begin{prop}
\label{LAPtestsuralpha} Let the assumptions of
Proposition \ref{LANproperty} and Theorem~\ref
{theoconsistenceJensenRahbek} be satisfied.
For testing (\ref{testingasymmetry}), the test defined by the
rejection region (\ref{test3}) has the local asymptotic power
\begin{eqnarray*}
\lim_{n \to
\infty}P_{n,\btau}^{\mathrm{S}}
\bigl({\mathrm C}^{\mathrm
{S}} \bigr)&=&1-\Phi \biggl\{\Phi^{-1}
\biggl(1-\frac{\underline
{\alpha}}{2} \biggr)-\frac{\tau_2-\tau_3}{\sigma_{T^\mathrm
{S}}} \biggr\}
\\
&&{} + \Phi \biggl\{-\Phi^{-1} \biggl(\frac{\underline{\alpha}}{2} \biggr)-
\frac{\tau_2-\tau_3}{\sigma_{T^\mathrm{S}}} \biggr\},
\end{eqnarray*}
where, recalling the notation $\be'=(1,-1,0)$,
\[
\sigma^2_{T^\mathrm{S}}= \cases{ ({\kappa}_{\eta}-1)
\be'{\cal I}_*^{-1}\be, &\quad when $\gamma_0< 0$,
\cr
({\kappa}_{\eta}-1)\be'{\cal I}^{-1}\be, &\quad
when $\gamma_0\geq0$.}
\]
\end{prop}

\subsection{Optimality issues}
We discuss, in this section, the optimality of the symmetry test
defined in
(\ref{test3}).
Let $\theta_0=(\omega_0,\alpha_0,\alpha_0,\beta_0)'$ be a parameter
value corresponding to a symmetric GARCH. Assume that, at this
point, $\gamma_0\geq0$. If $\gamma_0< 0$, it suffices to replace
${\cal I}$ by ${\cal I}_*$ in the sequel. A sequence of local
alternatives to this symmetric parameter is defined by
$\theta_0+\btau/\sqrt{n}$ where
$\btau'=(\tau_1,\tau_2,\tau_3,\tau_4)'$ is such that $\tau_2\neq
\tau_3$. Relations (\ref{lan})--(\ref{Infomat}) imply that
\[
\Lambda_{n,f}(\theta_0+\btau/\sqrt{n},
\theta_0)\stackrel {d} {\longrightarrow} {\cal N} \biggl(-
\frac{\iota_f}{8}\tilde{\btau}'{\cal I}\tilde{\btau},
\frac{\iota_f}{4}\tilde{\btau}'{\cal I}\tilde{\btau} \biggr) \qquad\mbox{under }P_{\theta_0}
\]
with
$\tilde{\btau}=(\tau_2,\tau_3,\tau_4)'$, which is the distribution
of the log-likelihood ratio in the statistical model ${\cal
N} \{\tilde{\btau},4{\cal I}^{-1}/\iota_{f} \}$ of
parameter $\tilde{\btau}$. In other words, denoting by $\tilde{\cal
T}$ a subset of $\mathbb{R}^3$ containing a neighborhood of
$\bzero$, for any $\tau_1$, the so-called local experiments
$ \{L_{n,f}(\theta_0+(\tau_1,\tilde{\btau}')/\sqrt{n}),\tilde
{\btau}\in
\tilde{\cal T} \}$ converge to the Gaussian experiment
$ \{{\cal N} (\tilde{\btau},4{\cal
I}^{-1}/\iota_{f} ),\tilde{\btau}\in\tilde{\cal T} \}$.

The asymmetry test (\ref{testingasymmetry})
corresponds
to the test
\[
\be'\tilde{\btau}= 0 \quad\mbox{against}\quad \be'\tilde{
\btau } \neq0
\]
in the limiting experiment.
The uniformly most powerful unbiased (UMPU) test based on
$\bX\sim{\cal N} (\tilde{\btau},4{\cal I}^{-1}/\iota
_{f} )$
is the test of rejection region
\[
C= \bigl\{\bigl|\be'\bX\bigr|/\sqrt{4\be'{\cal
I}^{-1}\be/\iota_{f}}>\Phi ^{-1}(1-\underline{
\alpha}/2) \bigr\}.
\]
This UMPU test has the power
%
\begin{equation}\quad
\label{bornepower} P_{\be'\tilde{\btau}}(C)=1-\Phi \biggl\{\Phi^{-1}
\biggl(1-\frac
{\underline{\alpha}}{2} \biggr)-c_{\be'\tilde{\btau}} \biggr\}+ \Phi \biggl\{-
\Phi^{-1} \biggl(\frac{\underline{\alpha}}{2} \biggr)-c_{\be'\tilde{\btau}} \biggr\}
\end{equation}
with\vspace*{1pt}
$c_{\be'\tilde{\btau}}=\frac{\be'\tilde{\btau}\sqrt{\iota
_{f}}}{2\sqrt{\be'{\cal
I}^{-1}\be}}$.
A test of (\ref{testingasymmetry}) whose level converges to
$\underline{\alpha}$, which is asymptotically unbiased, and whose
power converges to the bound in (\ref{bornepower}) will be called
asymptotically locally UMPU.
%
\begin{prop}
\label{coro2} Under the assumptions of Proposition \ref{LAP2tests},
the test (\ref{test3}) is asymptotically locally UMPU for the
testing problem (\ref{testingasymmetry}) if and only if the density of
$\eta_t$ has the form
%
\begin{equation}
\label{casqmleefficient}\qquad f(y)=\frac{a^a}{\Gamma(a)}e^{-a
y^2}|y|^{2a-1},\qquad
a>0,\qquad \Gamma(a)=\int_0^{\infty}t^{a-1}e^{-t}\,dt.
\end{equation}
\end{prop}
A figure displaying the density (\ref{casqmleefficient}) for
different values of $a$ is in the supplementary file [\citet{FraZak}].
Note that the
Gaussian density is obtained for $a=1/2$. The result was expected
because the ${\mathrm C}^{\mathrm{S}}$-test is based on the QMLE of
$\theta_{0}$, and the QMLE is obviously efficient in the Gaussian
case.
It can be shown that when the distribution of $\eta_t$ is of the
form (\ref{casqmleefficient}), the MLE does not depend on $a$. The
QMLE is then equal to the MLE, which makes obvious the ``if part'' of
Proposition \ref{coro2}. The ``only if'' part of the proposition shows
that there is necessarily an efficiency loss when the test is not
based on the MLE of $\theta_0$.

%
\begin{figure}

\includegraphics{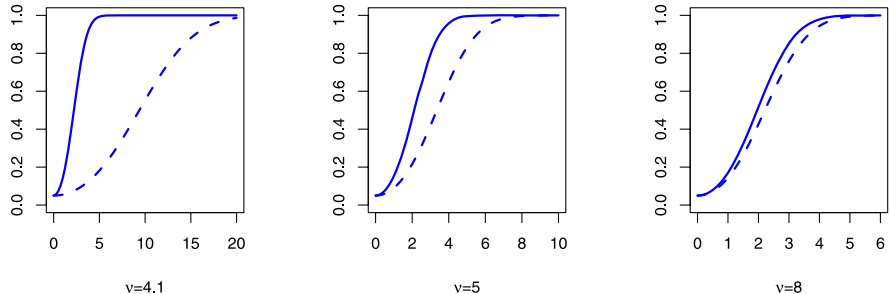}

\caption{Optimal asymptotic power (\protect\ref{bornepower}) (in
full line) and
local asymptotic power of the asymmetry test (\protect\ref{test3})
(in dotted
line) when $\eta_t$ follows a standardized Student distribution with
$\nu$ degrees of freedom. The horizontal axis correspond to the local
parameter $\be'\btau$.}\label{comparpower}
\end{figure}

This point is illustrated by Figure~\ref{comparpower}, in which the
local asymptotic power of the asymmetry test (in dotted lines) is
compared to the optimal asymptotic power given by
(\ref{bornepower}). In this figure, the noise $\eta_t$ is assumed to
satisfy a Student distribution with $\nu>2$ degrees of freedom,
standardized in such a way that $E\eta_t^2=1$. The parameters of the
model under the null are $\alpha_{0+}=\alpha_{0-}=0.2$,
$\beta_0=0.9$ and $\delta=1$, which corresponds to a nonstationary
model with $\gamma_0=0.045$. In the figure, it can be seen that the
local asymptotic power is far from the optimal power when $\nu$ is
small, but, as expected, the discrepancy decreases as $\nu$
increases.

\section{\texorpdfstring{Estimation when the power $\delta$ is unknown}{Estimation when the power delta is unknown}}
\label{Sec4} In this section, we consider the case where the power
$\delta$, now denoted $\delta_0$, is unknown and is jointly
estimated with $\theta_0$.
We rewrite the vector of parameters as $\zeta:= (\delta,\theta')'$,
which is assumed to belong to a compact parameter space $\Upsilon
\subset(0,\infty)^2 \times[0,\infty)^3$. The true parameters value
is denoted by
$\zeta_0:= (\delta_0,\theta_0')'$.
A QMLE of $\zeta$ is defined as any measurable solution
$\hat{\zeta}_n$ of
%
\begin{equation}
\label{qmldlta}\quad \hat{\zeta}_n=\bigl(\hat{\delta}_n, \hat{
\theta}_n'\bigr)'=\argmin_{\zeta
\in\Upsilon}
\frac{1}{n}\sum_{t=1}^n
\ell_t(\zeta),\qquad \ell_t(\zeta)=\frac{\epsilon_t^2}{\sigma_t^2(\zeta)}+\log\sigma
_t^2(\zeta),
\end{equation}
where
%
\begin{equation}
\label{varsigma} \sigma_t={\sigma}_t(\zeta)= \bigl(
\omega+\alpha_{+}\bigl(\epsilon _{t-1}^+\bigr)^{\delta}
+\alpha_{-}\bigl(-\epsilon_{t-1}^-\bigr)^{\delta} +
\beta{\sigma}_{t-1}^\delta(\zeta) \bigr)^{1/\delta}
\end{equation}
for
$t=1,\ldots,n$ [with initial values for $\epsilon_0$ and
$\sigma_0(\zeta)$]. The rescaled residuals are defined by
$\hat{\eta}_t=\eta_t(\hat{\zeta}_n)$ where
$\eta_t(\zeta)=\epsilon_t/\sigma_t(\zeta)$ for $t=1,\ldots,n$.
For identifiability reasons, we need to slightly reinforce
assumption A1 as follows.

{A3:}
The support of $\eta_t$ contains at least three points of the same
sign, and at least two points of opposite signs.

We also introduce the following technical assumption to handle the
derivatives of $\ell_t$ with respect to the exponent $\delta$.

{A4:}
$\forall\zeta\in\Upsilon$, $\beta<\|1/a_0^2(\eta_1)\|_{p}^{-1}$
and $\||\eta_1|^{\delta}\log|\eta_1|\|_p<\infty$
for some $p>1$.

For brevity, we only present results for the nonstationary cases.
%
\begin{theo}
\label{theoconsistenceJensenRahbekWithDelta} Let
(\ref{ARCH1deJensenRahbek})--(\ref{eta}) and \textup{A3} hold. Then the
QMLE defined in (\ref{qmldlta}) satisfies the following properties:
%
\begin{enumerate}[(ii)]
%
\item[(i)] \textup{Explosive case.}
When $\gamma_0>0$, if $P(\eta_1=0)=0$
\[
\bigl(\delta_n, \hat{\vartheta}_n'\bigr)
\to \bigl(\delta_0,\vartheta_0'\bigr)\qquad
\mbox{a.s. as } n\to\infty.
\]
If, in addition,
$\kappa_{\eta}\in(1,\infty)$,
$E|\log\eta_1^2|<\infty$, $\zeta_0\in
\mathring{\Upsilon}$, and \textup{A4} holds, then
%
\begin{equation}
\label{normalitydeleQMVdealphabis} \sqrt{n} \bigl(\bigl(\hat{
\delta}_n,\hat{\vartheta}_n'\bigr)-\bigl(
\delta _0,\vartheta_0'\bigr)
\bigr)'\stackrel{d} {\to} {\cal N} \bigl\{0,(\kappa
_{\eta}-1){\cal I}_{\delta}^{-1} \bigr\} 
\end{equation}
as $n\to\infty$, where ${\cal I}_{\delta}$ is a positive
definite matrix (see Lemma 3.1).
\item[(ii)] \textup{At the boundary of the stationarity region.}
When $\gamma_0= 0$, if $P(\eta_1=0)=0$, and
$\forall\zeta\in\Upsilon$, $\beta<\|1/a_0(\eta_1)\|_{p}^{-1}$
for some $p>1$,
%
\[
\bigl(\delta_n, \hat{\vartheta}_n'\bigr)
\to \bigl(\delta_0,\vartheta_0'\bigr)
\qquad\mbox{in probability as } n\to\infty.
\]
%
If, in addition,
$\zeta_0\in \mathring{\Upsilon}$,
$\kappa_{\eta}\in(1,\infty)$, $E|\log\eta_1^2|<\infty$
and \textup{A2} and \textup{A4} are satisfied,
then (\ref{normalitydeleQMVdealphabis}) holds.
\end{enumerate}
\end{theo}
The presence of parameter $\delta$ induces specific difficulties. It
turns out that the derivative of the criterion with respect to
$\delta$ involves the process
$({\partial\sigma_t^\delta}/{\partial\delta}-\log\sigma_t)$. A
strictly stationary approximation to this process can then be
obtained, but in a more complicated way than for the other
parameters. To save space, the proofs of this section are given in
the supplementary file [\citet{FraZak}].

Obviously, stationarity and symmetry tests could be derived as in
Sections~\ref{Sectest} and~\ref{secLAN}. Other tests concerning the
exponent $\delta$ [e.g., testing the TARCH model \mbox{($\delta=1$)}
against the GJR model ($\delta=2$)] could be considered as well, but
we leave this for further investigation.

\section{Proofs and complementary results}\label{complement}

\mbox{}

\begin{pf*}{Proof of Proposition \ref{convtoinfty}}
Writing $\omega_t=\omega(\xi_t)$ and $a_t=a(\xi_t)$, we have, for
all $t>1$ and
$1\leq k<t$,
%
\begin{equation}
\label{recurencesurh} h_t = \omega_{t-1}+\sum
_{j=1}^k\omega_{t-j-1}\prod
_{i=1}^ja_{t-i}+h_{t-k-1}\prod
_{i=1}^{k+1}a_{t-i}.
\end{equation}
We begin by showing (i). Since all the random variables involved
in (\ref{recurencesurh}) are positive, $h_t\geq
\underline{\omega}\prod_{i=1}^{t-1}a_{t-i}$. For any constant
$\rho> e^{-\gamma}$, we thus have, a.s.
\[
\liminf_{t\to\infty}\frac{1}{t}\log\rho^th_t
\geq\log\rho+ \lim_{t\to\infty}\frac{1}{t} \Biggl\{\log
\underline{\omega}+\sum_{i=1}^{t-1}\log
a_i \Biggr\} =\log\rho+\gamma>0
\]
by the ergodic theorem. It follows that $\log\rho^th_t$, and hence
$\rho^th_t$, tend to $+\infty$ a.s. as $n\to\infty$. The second
convergence is shown in just the same way, arguing that $E|\log
\xi_1^2|<\infty$ entails $\log\xi_t^2/t \to0$ a.s. as $t\to\infty$.

To show (ii), first consider the case where $h_0=0$. Note
that, for all $t$, the distribution of
$h_t=h_t(\xi_0,\ldots,\xi_{t-1})$ is equal to that of
%
\begin{equation}
h^*_t:=h_t(\xi_{t},\ldots,\xi_1)=
\omega_{1}+\sum_{j=1}^{t-1}
\omega_{j+1}\prod_{i=1}^ja_{i}.
\end{equation}
Note that, contrary to $(h_t)$, the sequence $(h^*_t)$ increases
with $t$. The Chung--Fuchs theorem applied to the random walk
$\sum_{i=1}^t\log a_{i}$ entails that\break $\limsup_{t\to\infty}
\prod_{i=1}^ta_{i}=+\infty$ a.s. It follows that $h^*_t\to+\infty$
as $t\to\infty$. We thus have $P(h_t\geq A)= P(h_t^*\geq A)\to1$
for all $A>0$, from which the first part of (ii) easily
follows.
To prove the first convergence of (\ref{empmeans}), note that the
dominated convergence theorem entails
\[
E\psi(h_t)=\int_0^{\infty}P \bigl
\{h_t^*<\psi^{-1}(u) \bigr\}\,du \to\int_0^{\infty}
\lim_{t\to\infty}
P \bigl\{h_t^*<
\psi^{-1}(u) \bigr\}\,du=0.
\]
The second convergence is shown similarly.
Now consider the case where the initial value is not equal to zero.
It is clear from (\ref{recurencesurh}), with $k=t-1$, that $h_t$ is
an increasing function of $h_0$. So the convergences to infinity
obtained when $h_0=0$, and the convergences in
(\ref{empmeans}), hold a fortiori when $h_0>0$.
\end{pf*}



\subsection{\texorpdfstring{Asymptotic behavior of the QMLE of $\vartheta_0$}{Asymptotic behavior of the QMLE of theta 0}}\label{SecA2}


Define the
$[0,\infty]$-valued process
\[
v_{t}(\vartheta)=\sum_{j=1}^{\infty}
\frac{\{\alpha_+(\eta
_{t-j}^+)^{\delta}+\alpha_-(-\eta_{t-j}^-)^{\delta}\}}{a_0(\eta
_{t-j})}\prod_{k=1}^{j-1}
\frac{\beta}{a_0(\eta_{t-k})} 
\]
with the convention $\prod_{k=1}^{j-1}=1$ when $j\leq1$. 
Let $\Theta_0=\{\theta\in\Theta\dvtx  \beta<e^{\gamma_0}\}$ and
$\Theta_p=\{\theta\in[0,\infty)^4\dvtx
\beta<\|1/a_0(\eta_1)\|_{p}^{-1}\}$.
%
\begin{lem}
\label{lemderivbeta}
\textup{(i)} When $\gamma_0>0$, for any $\theta\in\Theta_0$ the
process $v_{t}(\vartheta)$ is stationary and ergodic. Moreover, for
any compact $\Theta_0^*\subset\Theta_0$,
\[
\sup_{\theta\in\Theta_0^*} \biggl\llvert \frac{\sigma_t^{\delta}(\theta)}{h_t}-v_{t}(
\vartheta )\biggr\rrvert \to0\qquad\mbox{a.s. as $t\to\infty$}.
\]
Finally, for any $\theta\notin\Theta_0$ it holds that
${\sigma_t^{\delta}(\theta)}/{h_t}\to\infty$ a.s.

\textup{(ii)} When $\gamma_0=0$, for any $\theta\in\Theta_p$ with
$p\geq1$, the process $v_{t}(\vartheta)$ is stationary and
ergodic. Moreover, for any compact $\Theta_p^*\subset\Theta_p$,
\[
\sup_{\theta\in\Theta_p^*}\biggl\llvert \frac{\sigma_t^{\delta}(\theta
)}{h_t}-v_{t}(
\vartheta)\biggr\rrvert \to0\qquad\mbox{in $L^p$}.
\]
\end{lem}
\begin{pf}
Assuming, with no generality loss, that
$\sigma_0(\theta)=0$, we have
$\sigma_t^{\delta}(\theta)=\sum_{j=1}^t\beta^{j-1}z_{t-j}$ where
$z_t=\omega+\alpha_+(\epsilon_{t}^+)^{\delta}+\alpha_-(-\epsilon
_{t}^-)^{\delta}$
and
%
\begin{equation}\label{aetb}
\frac{\sigma_t^{\delta}(\theta)}{h_t}=\sum_{j=1}^t\beta
^{j-1} \Biggl\{\prod_{k=1}^j
\frac{h_{t-k}}{h_{t-k+1}} \Biggr\} \frac{z_{t-j}}{h_{t-j}}.
\end{equation}
Noting that
%
\begin{equation}
\label{h2surh1} \frac{h_{t-k}}{h_{t-k+1}}=\frac{h_{t-k}}{\omega_0+a_0(\eta
_{t-k})h_{t-k}}\leq\frac{1}{a_0(\eta_{t-k})},
\end{equation}
the rest of the proof follows from arguments similar to those used
in the proof of Lemma A.1 in FZ. Therefore is it omitted.
\end{pf}

\begin{lem}
\label{lemidentif} If 
$\theta\in\Theta_0$, we have $ v_{t}(\vartheta)= 1$,
a.s. if and only if $\vartheta=\vartheta_{0}$.
\end{lem}
\begin{pf}
Straightforward algebra shows that
%
\begin{equation}
\label{biensympa} v_{t}(\vartheta)a_0(
\eta_{t-1})= 
\beta v_{t-1}(\vartheta)+
\alpha_+\bigl(\eta_{t-1}^+\bigr)^{\delta}+\alpha _-\bigl(-
\eta_{t-1}^-\bigr)^{\delta}.
\end{equation}
Hence
\begin{eqnarray*}
\bigl\{v_{t}(\vartheta)-1\bigr\}a_0(
\eta_{t-1})&=&\beta v_{t-1}(\vartheta)-\beta_0+ (
\alpha_+-\alpha_{0+}) \bigl(\eta_{t-1}^+\bigr)^{\delta}\\
&&{}+(
\alpha_--\alpha _{0-}) \bigl(-\eta_{t-1}^-\bigr)^{\delta}.
\end{eqnarray*}
It follows that $v_{t}(\vartheta)= 1$ a.s. if and only if
\[
\beta -\beta_0+(\alpha_+-\alpha_{0+}) \bigl(
\eta_{t-1}^+\bigr)^{\delta}+(\alpha _--\alpha_{0-})
\bigl(-\eta_{t-1}^-\bigr)^{\delta}=0.
\]
%
Thus if $\vartheta\ne\vartheta_0$, $\eta_t$ takes at most two
values of different signs, in contradiction with assumption
A1. The conclusion follows.
\end{pf}

Let $\underline{\omega}=\inf\{\omega\mid\theta\in\Theta\}$,
$\underline{\alpha}=\inf\{\alpha_+, \alpha_-\mid\theta\in
\Theta\}$, $\underline{\beta}=\inf\{\beta\mid\theta\in\Theta\}$,
$\overline{\omega}=\sup\{\omega\mid\theta\in\Theta\}$,
$\overline{\alpha}=\sup\{\alpha_+, \alpha_-\mid\theta\in
\Theta\}$, $\overline{\beta}=\sup\{\beta\mid\theta\in\Theta\}$.
Denote\vspace*{2pt} by $K$ any constant whose value is unimportant and can change
throughout the proofs. Let $\check{\Theta}$ be the compact set of
the $\vartheta$'s such that $(\omega,\vartheta')'\in\Theta$.
%
\begin{lem}\label{lemv}
Suppose that $P(\eta_t=0)=0$. Then, for any $k>0$,
\[
E\sup_{\vartheta\in\check{\Theta}} \biggl(\frac1{v_t(\vartheta )}
\biggr)^k<\infty \quad\mbox{and}\quad E \sup_{\theta\in\Theta} \biggl(
\frac{h_t}{\sigma_t^{\delta
}(\theta)} \biggr)^k<\infty.
\]
\end{lem}
\begin{pf}
Let $\varepsilon>0$ such that
$p(\varepsilon):=P(|\eta_{t}|\leq\varepsilon)\in[0,1)$. If
$|\eta_{t-1}|>\varepsilon$, since the sum $v_{t}(\vartheta)$ is
greater than its first term, we have
%
\[
\frac{1}{v_{t}(\vartheta)}\leq\frac{a_0(\eta_{t-1})}{\alpha
_+(\eta_{t-1}^+)^{\delta}+\alpha_-(-\eta_{t-1}^-)^{\delta}} \leq
\frac{\max(\alpha_{0+},\alpha_{0-}) }{\underline{\alpha
}}+\frac{\beta_0}{\underline{\alpha}\varepsilon^{\delta
}}:=K(\varepsilon).
\]
%
Iterating this method, we can write
\[
\sup_{\vartheta\in\check{\Theta}} \frac{1}{v_{t}(\vartheta)}\leq K(\varepsilon) \sum
_{i=1}^{\infty}\mathbh1_{|\eta_{t-1}|\leq\varepsilon}\cdots
\mathbh1_{|\eta_{t-i+1}|\leq\varepsilon} \mathbh1_{|\eta_{t-i}|> \varepsilon} 
\biggl(
\frac{\overline{a}_0(\varepsilon)}{\underline{\beta}} \biggr)^{i-1},
\]
where
$\overline{a}_0(\varepsilon)=\max(\alpha_{0+},\alpha_{0-})\epsilon
^{\delta}+\beta_0$.
It follows that, for any integer $k$, 
\[
E\sup_{\vartheta\in\check{\Theta}} \biggl(\frac1{v_t( \vartheta )}
\biggr)^k\leq \bigl\{K(\varepsilon)\bigr
\}^k\bigl\{1-p(\varepsilon)\bigr\} \sum_{i=1}^{\infty}
p(\varepsilon)^{i-1}
\biggl(\frac{\overline{a}_0(\varepsilon)}{\underline{\beta}}
\biggr)^{k(i-1)}.
\]
Noting that $\lim_{\varepsilon\to0}p(\varepsilon)=0$ and
$\lim_{\varepsilon\to0}\overline{a}_0(\varepsilon)=\beta_0$, we have
$p(\varepsilon)
(\frac{\overline{a}_0(\varepsilon)}{\underline{\beta}} )^{k}<1$
for $\varepsilon$ sufficiently small. The first result of the lemma
is thus proven.

Similarly, we have for $|\eta_{t-1}|>\varepsilon$,
\[
\frac{h_t}{\sigma_{t}^{\delta}(\theta)} 
\leq
\frac{\omega_0}{\underline{\omega}}+\frac{\overline
{\alpha}}{\underline{\alpha}}+\frac{\beta_0}{\underline{\alpha
}\varepsilon^{\delta}}:=H(\varepsilon)
\]
%
and for $|\eta_{t-1}|\leq\varepsilon$ and
$|\eta_{t-2}|>\varepsilon$,
\[
\frac{h_t}{\sigma_{t}^{\delta}(\theta)} \leq\frac{\omega_0}{\underline{\omega}}+\frac{\overline
{a}_0(\varepsilon)}{\underline{\beta}}H(\varepsilon).
\]
More generally,
\begin{eqnarray*}
\sup_{\theta\in\Theta}\frac{h_t}{\sigma_{t}^{\delta}(\theta
)}&\leq&
\sum
_{i=1}^{\infty}\mathbh1_{|\eta_{t-1}|\leq\varepsilon}\cdots
\mathbh1_{|\eta_{t-i+1}|\leq\varepsilon} \mathbh1_{|\eta_{t-i}|> \varepsilon}
\\
&&\hspace*{12.9pt}{} \times 
\Biggl(\frac{\omega_0}{\underline{\omega}}\sum
_{j=0}^{i-2} \biggl(\frac{\overline{a}_0(\varepsilon)}{\underline{\beta}}
\biggr)^j+ \biggl(\frac{\overline{a}_0(\varepsilon)}{\underline{\beta}} \biggr)^{i-1}H(
\varepsilon) \Biggr).
\end{eqnarray*}
The conclusion follows by the same arguments as before.
\end{pf}

\begin{pf*}{Proof of the consistency results in cases \textup{(ii)} and
\textup{(iii)} of Theorem~\ref{theoconsistenceJensenRahbek}} Note that $(\hat{\omega}_n,\hat
{\vartheta}_n')=
\argmin_{\theta\in\Theta}Q_n(\theta)$, where
$Q_n(\theta)=n^{-1}\*\sum_{t=1}^n \{\ell_t(\theta)-\ell
_t(\theta_0) \}$.
We have
\[
Q_n(\theta)=\frac{1}{n}\sum_{t=1}^n
\eta_t^2 \biggl\{ \biggl(\frac
{h_t}{\sigma_t^{\delta}(\theta)}
\biggr)^{2/{\delta}}-1 \biggr\}+ \log \biggl(\frac{\sigma_t^{\delta}(\theta)}{h_t}
\biggr)^{2/{\delta}} =O_n(\vartheta)+ R_n(\theta),
\]
where
\[
O_n(\vartheta)= \frac{1}{n}\sum_{t=1}^n
\eta_t^2 \biggl\{\frac{1}{v_{t}^{2/{\delta
}}(\vartheta)}-1 \biggr\}+\log
v_{t}^{2/{\delta}}(\vartheta)
\]
and
\[
R_n(\theta)=\frac{1}{n}\sum_{t=1}^n
\eta_t^2 \biggl\{ \biggl(\frac
{h_t}{\sigma_t^{\delta}(\theta)}
\biggr)^{2/{\delta}}- \frac{1}{v_{t}^{2/{\delta}}(\vartheta)} \biggr\}+\log \biggl(\frac
{\sigma_t^{\delta}(\theta)}{h_tv_{t}(\vartheta)}
\biggr)^{2/{\delta}}.
\]
It suffices to consider the case $\theta\in\Theta^*_0$ where
$\Theta^*_0$ is an arbitrary compact subset of~$\Theta_0$, because
by Lemma \ref{lemderivbeta}(i) $Q_n(\theta)\to\infty$ a.s. if
$\theta\notin\Theta_0$.
We have by stationarity and ergodicity
of $v_{t}(\vartheta)$, a.s.
\[
\lim_{n\to\infty}O_n(\vartheta)= E \biggl\{
\frac{1}{v_{1}^{2/{\delta}}(\vartheta)}-1+\log v_{1}^{2/{\delta}}(\vartheta)
\biggr\}\geq0,
\]
because $\log x \leq x-1$ for $x>0$. The inequality is strict except
when $v_{1}(\vartheta)=1$ a.s. By Lemma \ref{lemidentif} we thus
have $E\{O_n(\vartheta)\}\geq0$, with equality only if
$\vartheta=\vartheta_0$.

By Lemma \ref{lemv} we prove, as in FZ, that
%
\begin{equation}
\label{res1preuvetheoconsistenceJensenRahbekII}\qquad \lim_{n\to\infty}\sup
_{\theta\in\Theta_0^*}\bigl\llvert R_n(\theta )\bigr\rrvert = 0
\qquad\mbox{a.s. }\Bigl[\mbox{resp., } \lim_{n\to\infty}\sup_{\theta\in\Theta_p^*}
\bigl\llvert R_n(\theta )\bigr\rrvert = 0 \mbox{ in } L^1\Bigr],
\end{equation}
when $\gamma_0>0$ (resp., $\gamma_0=0$) and $\Theta_0^*, \Theta_p^*$
are defined in Lemma \ref{lemderivbeta}, which completes the proof.
%
\end{pf*}

We now need to introduce new $[0,\infty]$-valued processes. Let
$a(\eta_{t})=\break\alpha_+(\eta_t^+)^{\delta}+\alpha_-(-\eta
_t^-)^{\delta}+\beta$
and
\begin{eqnarray*}
d_{t}^{\alpha_+}&=&\sum_{j=1}^{\infty}
\frac{(\eta_{t-j}^+)^{\delta
}}{a_0(\eta_{t-j})}\prod_{k=1}^{j-1}
\frac{\beta_0}{a_0(\eta
_{t-k})},\qquad d_{t}^{\alpha_-}=\sum
_{j=1}^{\infty}\frac{(-\eta_{t-j}^-)^{\delta
}}{a_0(\eta_{t-j})}\prod
_{k=1}^{j-1}\frac{\beta_0}{a_0(\eta_{t-k})},
\\
d_{t}^{\beta}&=&\sum_{j=2}^{\infty}
\frac{(j-1)\{\alpha_{0+}(\eta
_{t-j}^+)^{\delta}+\alpha_{0-}(-\eta_{t-j}^-)^{\delta}\}}{\beta_0
a_0(\eta_{t-j})}\prod_{k=1}^{j-1}
\frac{\beta_0}{a_0(\eta_{t-k})}.
\end{eqnarray*}

\begin{lem}
\label{lemderivalphabeta} Assume $\gamma_0\geq0$ and
$E\eta_t^4<\infty$. We have
\[
\frac{1}{\sqrt{n}}\sum_{t=1}^n
\frac{\partial\ell_t}{\partial
\vartheta}(\theta_0) 
\stackrel{d} {\to} {\cal N} \bigl\{0,(\kappa_{\eta}-1){\cal I} \bigr\}
\qquad\mbox{as $n\to \infty$},
\]
where 
${\cal I}=\frac4{\delta^2}Ed_1d_1'$
and
$d_t'= (d_{t}^{\alpha_+}, d_{t}^{\alpha_-},
d_{t}^{\beta} )$.
Moreover, ${\cal I}$ is nonsingular.
\end{lem}
\begin{pf}
Since $E\log\beta_0/a_0(\eta_1)<0$, by the Cauchy root test,
the processes $d_{t}^{\alpha_+}, d_{t}^{\alpha_-}$ and $d_{t}^{\beta
}$ are stationary and ergodic.
Still assuming $\sigma_0^2=0$, we have
\begin{eqnarray*}
\frac{\partial\sigma_t^{\delta}}{\partial(\alpha_+,\alpha
_-)}(\theta)&=& \sum_{j=1}^t
\beta^{j-1}\bigl(\bigl\{\epsilon_{t-j}^+\bigr\}^{\delta},
\bigl\{-\epsilon _{t-j}^-\bigr\}^{\delta}\bigr),\\
\frac{\partial\sigma_t^2}{\partial\beta}(\theta)&=& \sum_{j=2}^t(j-1)
\beta^{j-2}z_{t-j}.
\end{eqnarray*}
Thus, using a direct
extension of (\ref{h2surh1}),
\begin{eqnarray*}
\frac{1}{\sigma^{\delta}_t(\theta_0)}\,\frac{\partial\sigma
_t^{\delta}}{\partial(\alpha_+,\alpha_-)}(\theta_0)&=& \sum
_{j=1}^t\beta^{j-1} \Biggl\{\prod
_{k=1}^j\frac{\sigma^{\delta
}_{t-k}(\theta_0)}{\sigma^{\delta}_{t-k+1}(\theta_0)} \Biggr\}
\frac{\{(\epsilon_{t-j}^+)^{\delta},(-\epsilon_{t-j}^-)^{\delta}\}
}{\sigma^{\delta}_{t-j}(\theta_0)}
\\
&\leq& \bigl(d_{t}^{\alpha_+}(\vartheta_0),
d_{t}^{\alpha_-}(\vartheta _0)\bigr),
\\
\frac{1}{\sigma^{\delta}_t(\theta_0)}\,\frac{\partial\sigma
_t^{\delta}}{\partial\beta}(\theta_0)&=& \sum
_{j=2}^t(j-1)\beta_0^{j-2}
\Biggl\{\prod_{k=1}^j\frac{\sigma
^{\delta}_{t-k}(\theta_0)}{\sigma^{\delta}_{t-k+1}(\theta
_0)}
\Biggr\}\frac{z_{t-j}}{\sigma^{\delta}_{t-j}(\theta_0)}
\\
&\leq &d_{t}^{\beta}(\vartheta_0),
\end{eqnarray*}
where the first inequality stands componentwise. Moreover, we have
\[
0\leq d_{t}^{\alpha_+}(\vartheta_0)-
\frac{1}{\sigma^{\delta
}_t}\frac{\partial\sigma_t^{\delta}}{\partial\alpha^+}(\theta _0)\leq
s_{t_0}+r_{t_0},
\]
where
\begin{eqnarray*}
s_{t_0}&=&\sum_{j=1}^{t_0}
\frac{(\eta_{t-j}^+)^{\delta}}{a_0(\eta
_{t-j})}\prod_{k=1}^{j-1}
\frac{\beta_0}{a_0(\eta_{t-k})}- \frac{(\epsilon_{t-j}^+)^{\delta}}{\beta_0\sigma^{\delta
}_{t-j}(\theta_0)}\prod_{k=1}^j
\frac{\beta_0\sigma^{\delta
}_{t-k}(\theta_0)}{\sigma^{\delta}_{t-k+1}(\theta_0)},
\\
r_{t_0}&=&\sum_{j=t_0+1}^{\infty}
\frac{(\eta_{t-j}^+)^{\delta
}}{a_0(\eta_{t-j})}\prod_{k=1}^{j-1}
\frac{\beta_0}{a_0(\eta_{t-k})}.
\end{eqnarray*}
For all $p\geq1$, $\|r_{t_0}\|_p\to0$ as $t_0\to\infty$ because
$\|\beta_0/a_0(\eta_1)\|_p<1$ and
$\|(\eta_1^+)^{\delta}/\break a_0(\eta_1)\|_p<1/\alpha_{0+}$. Since, in
addition, $\|\beta_0
\sigma_{t-1}^{\delta}(\theta_0)/\sigma^{\delta}_t(\theta_0)\|_p<1$,
and
\[
\biggl\llVert \frac{\beta_0}{a_0(\eta_{t-1})}-\frac{\beta_0 \sigma
_{t-1}^{\delta}(\theta_0)}{\sigma^{\delta}_t(\theta_0)}\biggr\rrVert _p=
\biggl\llVert \frac{\beta_0\omega_0}{a_0(\eta_{t-1})\sigma^{\delta
}_t(\theta_0)}\biggr\rrVert _p\to 0
\]
as $t\to\infty$ by the dominated convergence theorem,
$s_{t_0}=s_{t_0}(t)$ converges to 0 in $L^p$ as
$t\to\infty$. The same derivations hold true when $d_{t}^{\alpha_+}$
is replaced by $d_{t}^{\alpha_-}$ and~$d_{t}^{\beta}$. Therefore,
$d_{t}^{\alpha_+}, d_{t}^{\alpha_-}$ and $d_{t}^{\beta}$ have
moments of any order, and
%
\begin{equation}\label{derminusdt}
\biggl\llVert \frac{1}{\sigma^{\delta}_t}\,\frac{\partial
\sigma_t^{\delta}}{\partial\vartheta}(\theta_0)-
d_t\biggr\rrVert \to0
\end{equation}
in $L^p$ for any $p\geq1$.

Using (\ref{derminusdt}) and the ergodic theorem, we thus have, as
$n\to\infty$,
\[
\operatorname{Var}\frac{1}{\sqrt{n}}\sum_{t=1}^n
\frac{\partial}{\partial\vartheta}\ell_t(\theta_0)
=\frac{4}{\delta^{2}}
\frac{\kappa_\eta-1}{n}\sum_{t=1}^nE
\bigl(d_td_t'\bigr)+o(1) \to(
\kappa_{\eta}-1){\cal I}.
\]
Moreover, it can be shown as in FZ that the Lindeberg condition is
satisfied, allowing us to apply the Lindeberg central limit theorem for
martingale differences; see \citet{Bil95}, page 476.

Now we show that ${\cal I}$ is nonsingular. Suppose there
exists $x=(x_1,x_2,x_3)' \in\mathbb{R}^3$ such that $x'{\cal
I}x=0$. Then we get $x'd_t=0$, that is,
\begin{eqnarray*}
&&\sum_{j=1}^{\infty} \biggl(x_1
\frac{(\eta
_{t-j}^+)^{\delta}}{a(\eta_{t-j})}+ x_2\frac{(-\eta_{t-j}^-)^{\delta}}{a(\eta_{t-j})}+ x_3(j-1)
\frac{\alpha_+(\eta_{t-j}^+)^\delta+\alpha_-(-\eta
_{t-j}^-)^\delta}{\beta
a(\eta_{t-j})} \biggr)
\\
&&\qquad{} \times\prod_{k=1}^{j-1}
\frac
{\beta}{a(\eta_{t-k})}=0\qquad\mbox{a.s.}
\end{eqnarray*}
It follows that
$x_1(\eta_{t-1}^+)^{\delta}+
x_2(-\eta_{t-1}^-)^{\delta}=z_{t-2}$,
a.s. where $z_{t-2}$ is a measurable function of the $\eta_{t-j}$
with $j>1$. Because $\eta_{t-1}$ is independent of $z_{t-2}$, this
variable must be a.s. constant. In view of assumption A1, this
entails $x_1=x_2=0$ and then $x_3=0$. Therefore, ${\cal I}$ is
nonsingular.
\end{pf}


\begin{lem}
\label{lemintermediaresurladerivee} Let $\varpi$ be an arbitrary
compact subset of $[0,\infty)$. Assume that $E\log\eta_1^2<\infty$.
When $\gamma_0>0$ we have, a.s.
\begin{eqnarray*}
\sum_{t=1}^{\infty}\sup_{\theta\in\Theta_0}
\biggl\llvert \frac{\partial
}{\partial\omega}\ell_t(\theta)\biggr\rrvert &<&
\infty,\qquad
\sum_{t=1}^{\infty}
\sup_{\theta\in\Theta_0}\biggl\llVert \frac
{\partial^2}{\partial
\omega\,\partial
\theta}\ell_t(
\theta)\biggr\rrVert <\infty,
\\
\sup_{\omega\in\varpi} \Biggl\llvert \frac{1}{n}\sum
_{t=1}^n\frac{\partial^2 \ell_t(\omega,\vartheta_0)}{\partial
\theta_{i+1}\,\partial \theta_{j+1}}-{\cal I}_{ij} \Biggr\rrvert &=&o(1)
\qquad\mbox{for all }i,j\in \{1,2,3\},
\\
\frac{1}{n}\sum_{t=1}^n\sup
_{\theta\in\Theta} \biggl\llvert \frac{\partial^3}{\partial
\theta_i\,\partial
\theta_j\,\partial
\theta_k}\ell_t(
\theta)\biggr\rrvert &=&O(1) \qquad\mbox{for all }i,j,k\in \{2,3,4\}.
\end{eqnarray*}
When $\gamma_0=0$ we have, for all $ i,j,k\in\{2,3,4\}$,
%
\begin{eqnarray}
\label{traficdederivees5}
\sup_{\omega\in\varpi} \Biggl\llvert \frac{1}{n}\sum
_{t=1}^n\frac{\partial^2 \ell_t(\omega,\alpha_0,\beta_0)}{\partial
\theta_{i+1}\,\partial \theta_{j+1}}-{\cal I}_{ij} \Biggr\rrvert
&=&o_P(1),
\\
\label{traficdederivees6}
\frac{1}{n}\sum_{t=1}^n\sup
_{\theta\in\Theta_4} \biggl\llvert \frac{\partial^3}{\partial
\theta_i\,\partial
\theta_j\,\partial
\theta_k}\ell_t(
\theta)\biggr\rrvert &=&O_P(1).
\end{eqnarray}
\end{lem}
\begin{pf}
This is similar to that
of Lemma A.5. in FZ, therefore is it omitted.
\end{pf}

\begin{pf*}{Proof of the asymptotic normality in case \textup{(ii)} of
Theorem \ref{theoconsistenceJensenRahbek}}
An expansion of the criterion
derivative gives
%
\begin{equation}
\label{DLdanspreuvenormalitydeleQMVdealpha}
\pmatrix{ \frac{1}{\sqrt{n}}\sum
_{t=1}^n\frac{\partial
}{\partial
\omega}
\ell_t(\hat{\theta}_n)
\cr
0} = \frac{1}{\sqrt{n}}\sum
_{t=1}^n\frac{\partial
}{\partial
\theta}
\ell_t(\theta_0)+{\cal J}_n\sqrt{n}(\hat{
\theta}_n-\theta_0),
\end{equation}
where ${\cal J}_n$ is a $4\times4$ matrix whose elements have the
form
\[
{\cal J}_n(i,j)=\frac{1}{n}\sum_{t=1}^n
\frac{\partial^2}{\partial
\theta_i\,\partial\theta_j}\ell_t\bigl(\theta_{i}^*\bigr),
\]
where $\theta^*_{i}=(\omega_{i}^*,\alpha_{i+}^*,
\alpha_{i-}^*,\beta_{i}^*)'$ is between $\hat{\theta}_n$ and
$\theta_0$. Moreover, it can be shown that, for $i,j=1,2,3$,
%
\begin{equation}
\label{cvgdeJndanspreuvenormalitydeleQMVdealpha}\quad {\cal J}_n(i+1,1)=o(1/
\sqrt{n}),\qquad {\cal J}_n(i+1,j+1)\to{\cal I}(i,j) \qquad\mbox{a.s.}
\end{equation}
The conclusion follows from the last rows of
(\ref{DLdanspreuvenormalitydeleQMVdealpha}) and Lemma
\ref{lemderivalphabeta}.
\end{pf*}

\begin{pf*}{Proof of the asymptotic normality in case (iii) of
Theorem~\ref{theoconsistenceJensenRahbek}} 
Note that (\ref{DLdanspreuvenormalitydeleQMVdealpha}) continues to
hold. In view of (\ref{traficdederivees5})--(\ref
{traficdederivees6}), we have
\[
{\cal J}_n(i+1,j+1)\to{\cal I}(i,j) \qquad\mbox{in probability as } n\to
\infty.
\]
To conclude, by the arguments used in case (ii), it suffices to
show that
%
\begin{equation}
\label{convto0}\quad\mbox{for $i=2,3,4$}\qquad E\bigl|{\cal J}_n(i,1)\sqrt{n}(
\hat{\omega}_n-\omega_0)\bigr|\to0 \qquad\mbox{as }n\to\infty.
\end{equation}
Noting that
%
\begin{equation}
\label{maj20} \frac{1}{\sigma_t^{\delta}(\theta)}\sum_{j=1}^t
\beta ^{j-1}\bigl(\epsilon_{t-j}^+\bigr)^{\delta}\leq
\frac{1}{\alpha_+},
\end{equation}
and $\beta_2^*<1$ for $n$ large enough, and using the compactness of
$\Theta$, we
obtain
\begin{eqnarray*}
&&\bigl|{\cal J}_n(2,1)\sqrt{n}(\hat{\omega}_n-
\omega_0)\bigr|
\\
&&\qquad\leq \frac{K}{\sqrt{n}}\sum_{t=1}^n
\biggl(\frac{2h_t^{2/\delta}\eta_t^2}{\sigma_t^2(\theta
^*_2)}+1 \biggr) \frac{ \{\sum_{j=1}^t (\beta_2^* )^{j-1}(\epsilon
_{t-j}^+)^{\delta} \}
\{\sum_{j=1}^t (\beta_2^* )^{j-1} \}}{\sigma
_t^{2\delta}(\theta^*_2)}
\\
&&\qquad\leq \frac{K}{\sqrt{n}}\sum_{t=1}^n
\biggl(\frac{2h_t^{2/\delta}\eta_t^2}{\sigma_t^2(\theta
^*_2)}+1 \biggr)\frac{h_t}{\sigma_t^{\delta}(\theta^*_2)} \frac{1}{h_{t}}.
\end{eqnarray*}
Hence, by Lemma \ref{lemv} and H\"{o}lder's inequality
\[
E\bigl|{\cal J}_n(2,1)\sqrt{n}(\hat{\omega}_n-
\omega_0)\bigr| \leq \frac{K}{\sqrt{n}}\sum
_{t=1}^n E\frac{1}{h_{t}^{1+\tau}}
\]
for any $\tau>0$. The same
bound is obtained when ${\cal J}_n(2,1)$ is replaced by ${\cal
J}_n(3,1)$ and ${\cal J}_n(4,1)$. Moreover,
\[
h_t=\omega_0(1+Z_{t-1}+Z_{t-1}Z_{t-2}+
\cdots+ Z_{t-1}\cdots Z_{1})+ Z_{t-1}\cdots
Z_{0}\sigma_0^2.
\]
Hence
\[
\frac{1}{h_{t}^{1+\tau}}\leq\frac{1}{\omega_0^{1+\tau
}(1+Z_{t-1}+Z_{t-1}Z_{t-2}+ \cdots+ Z_{t-1}\cdots Z_{1})}.
\]
By assumption A2, the conclusion follows. 
\end{pf*}

\begin{pf*}{Proof of Theorem \ref{propo1}} To save space, this is
displayed in the supplementary file [\citet{FraZak}].
\end{pf*}

\subsection{Stationarity test}

\mbox{}

\begin{pf*}{Proof of Theorem \ref{jointdistribution}}
In the stationary case $\gamma_0<0$, standard arguments show that
%
\begin{equation}
\label{eq1} \hat{\gamma}_n={\gamma}_n(
\theta_0)+ \frac{\partial
{\gamma}_n(\theta_0)}{\partial
\theta'}(\hat{\theta}_n-
\theta_0)+o_P\bigl(n^{-1/2}\bigr)
\end{equation}
with
%
\begin{eqnarray}
\label{dergamma}\quad \frac{\partial{\gamma}_n(\theta_0)}{\partial\theta}&=& \frac{-1}{n}\sum
_{t=1}^n \frac1{a_{0}(
\eta_t)} \left[ \bigl\{ a_0(\eta_t)-
\beta_0 \bigr\}\frac{1}{h_t}\frac{\partial
\sigma^{\delta}_t(\theta_0)}{\partial\theta}- \pmatrix{ 0
\cr
\bigl(\eta^+_t\bigr)^\delta
\cr
\bigl(-\eta^-_t
\bigr)^\delta
\cr
1} \right]
\nonumber\\[-8pt]\\[-8pt]
&=&-\Psi+o_P(1),
\nonumber
\end{eqnarray}
where
$\Psi=(1-\nu_1)\Omega-a $ and $\Omega=E_{\infty} \frac
{1}{h_t}\,\frac{\partial
\sigma_t^\delta(\theta_0)}{\partial\theta}$.
Moreover the QMLE satisfies
%
\begin{equation}
\label{eq2} \sqrt{n}(\hat{\theta}_n-\theta_0)=-{\cal
J}^{-1}\frac{1}{\sqrt {n}}\sum_{t=1}^n
\bigl(1-\eta_t^2\bigr) \frac{2}{\delta h_t}\,
\frac{\partial
\sigma_t^\delta(\theta_0)}{\partial\theta}+o_P(1).
\end{equation}
In view of (\ref{eq1}), (\ref{dergamma}) and (\ref{eq2}), we have
\[
\sqrt{n}(\hat{\gamma}_n-\gamma_0) =
\frac1{\sqrt{n}}\sum_{t=1}^n
u_t+\Psi'{\cal J}^{-1}\frac{1}{\sqrt{n}}\sum
_{t=1}^n\bigl(1-\eta_t^2
\bigr) \frac{2}{\delta h_t}\,\frac{\partial
\sigma_t^\delta(\theta_0)}{\partial\theta}+o_P(1).
\]
Note that
\[
\operatorname{Cov} \Biggl(\frac1{\sqrt{n}}\sum_{t=1}^n
u_t, \Psi'{\cal J}^{-1}\frac{1}{\sqrt{n}}\sum
_{t=1}^n\bigl(1-\eta_t^2
\bigr) \frac{2}{\delta h_t}\,\frac{\partial
\sigma_t^\delta(\theta_0)}{\partial\theta} \Biggr) =\frac{2c}{\delta}
\Omega'{\cal J}^{-1}\Psi,
\]
where $c=\operatorname{Cov}(u_t, 1-\eta_t^2)$. The Slutsky lemma and the
central limit theorem for martingale differences thus entail
\[
\sqrt{n}(\hat{\gamma}_n-\gamma_0) \stackrel{d} {\to}
{\cal N} \biggl(0,\sigma^2_u+4\frac{c}{\delta}
\Omega'{\cal J}^{-1}\Psi +(\kappa_{\eta}-1)
\Psi'{\cal J}^{-1}\Psi \biggr).
\]
Now let $\overline{\theta}_0=(\omega_0, \alpha_{0+}, \alpha_{0-},
0)'$.
Noting that $\overline{\theta}_0'\,\partial
\sigma_t^\delta(\theta_0)/\partial\theta=h_t$ almost surely, we
have
\[
E \biggl\{\frac{1}{h_t}\,\frac{\partial\sigma_t^\delta(\theta
_0)}{\partial\theta} \biggl(1-\frac{1}{h_t}
\frac{\partial\sigma
_t^\delta(\theta_0)}{\partial\theta'}\overline{\theta}_0 \biggr) \biggr\}=0,
\]
which entails $\frac{\delta^2}{4}{\cal J}\overline{\theta}_0=\Omega$
and $\Omega'{\cal J}^{-1}\Omega=\frac{\delta^2}{4}$. 
It follows that
\[
\Omega'{\cal J}^{-1}\Psi=(1-\nu_1)
\frac{\delta^2}{4}-\frac{\delta
^2}{4}\overline{\theta}_0'a=
\frac{\delta^2}{4}(1-\nu_1-\alpha_{0+}\tilde{
\nu}_{1+}-\alpha _{0-}\tilde{\nu}_{1-})= 0.
\]
%
We also have $\Psi'{\cal J}^{-1}\Psi=a' {\cal J}^{-1}a-(1-\nu_1)^2$,
which completes the proof of the asymptotic distribution
(\ref{amontrerenappendix3}) in the case $\gamma_0<0$.

Now consider the case $\gamma_0\geq0$. Let 
$\theta_n^*$ be a sequence such that $\|\theta_n^*-\theta_0\|\leq
\|\hat{\theta}_n-\theta_0\|$. By Proposition \ref{convtoinfty}
(using assumption A2 when $\gamma_0=0$), we have
\[
\frac{1}{\sqrt{n}}\sum_{t=1}^n
\frac{1}{\sigma_t^\delta(\theta
_n^*)}\frac{\partial\sigma_t^\delta(\theta_n^*)}{\partial\omega
}=o(1)\qquad\mbox{a.s. (resp., in probability) as
$n\to\infty$,}
\]
when $\gamma_0>0$ (resp., when $\gamma_0=0$). It can be deduced that,
under the same conditions, $\sqrt{n}\,\frac{\partial^2 \gamma
_n(\theta_n^*)}{\partial\omega\,\partial\theta}=o(1)$, and
$\sqrt{n}(\hat{\theta}-\theta_0)'\,\frac{\partial^2 \gamma
_n(\theta_n^*)}{\partial\theta\,\partial\theta'}(\hat{\theta
}-\theta_0)=o(1)$,\vspace*{1pt}
which entails that (\ref{eq1}) still holds. 
The previous arguments show that (\ref{dergamma}) holds with
\[
\Omega=E\pmatrix{0
\cr
d_t^{\alpha_+}(\theta_0)
\cr
d_t^{\alpha_-}(\theta _0)
\cr
d_t^{\beta}(
\theta_0)} =\frac{1}{1-\nu_1}\pmatrix{ 0
\cr
\tilde{
\nu}_{1+}
\cr
\tilde{\nu}_{1-}
\cr
\nu_1/\beta}
\quad\mbox{and}\quad \Psi=\pmatrix{ 0
\cr
0
\cr
0}.
\]
The conclusion follows.
\end{pf*}

\subsection{Asymptotic local powers}

\mbox{}

\begin{pf*}{Proof of Proposition \ref{LANproperty}}
The LAN of GARCH models has already been established in the stationary
case; see \citet{DroKla97}, \citet{LeeTan05}. The
nonstationary case will be studied under more general assumptions in
the proof of Proposition \ref{LANpropertybis}.
\end{pf*}

\begin{pf*}{Proof of Proposition \ref{LANpropertybis}}
Let the functions
\[
g_1(y)=1+y\frac{f'}{f}(y) \quad\mbox{and}\quad g_2(y)=1+2y
\frac
{f'}{f}(y)+y^2 \biggl(\frac{f'}{f}
\biggr)'(y).
\]
Introduce also the notation
\[
\Delta_{1,t}(\theta)=\frac{1}{\sigma_t(\theta)}\,\frac{\partial^2
\sigma_t(\theta)}{\partial\theta\,\partial\theta'},\qquad \Delta
_{2,t}(\theta)=\frac{1}{\delta^2\sigma_t^{2\delta}(\theta)}\,\frac
{\partial
\sigma_t^\delta(\theta)}{\partial\theta}\,
\frac{\partial
\sigma_t^\delta(\theta)}{\partial\theta'}.
\]
A Taylor
expansion of $\theta_n\mapsto\Lambda_{n,f}(\theta_n,\theta_0)$
around $\theta_0$ yields
%
\begin{equation}
\label{res1proofLAN} \Lambda_{n,f}(\theta_n,
\theta_0)={\btau'}S_{n,f}(
\theta_0)- \tfrac{1}{2}\btau'
\mathfrak{I}_n\bigl(\theta_n^*\bigr)\btau+{\cal
R}_n,
\end{equation}
where $\theta_n^*$ is between $\theta_0$ and $\theta_n$,
%
\begin{eqnarray}
\label{Sn}\quad S_{n,f}(\theta_0)&=&\frac{-1}{\sqrt{n}}\sum
_{t=1}^ng_1(\eta_t)
\frac{1}{\delta h_t}\,\frac{\partial
\sigma_t^\delta(\theta_0)}{\partial\theta},
\nonumber\\[-8pt]\\[-8pt]
\mathfrak{I}_n(\theta)&=&\frac{1}{n}\sum
_{t=1}^ng_1 \biggl(\frac
{\epsilon_t}{\sigma_t(\theta)}
\biggr)\Delta_{1,t}(\theta)- \frac{1}{n}\sum
_{t=1}^ng_2 \biggl(\frac{\epsilon_t}{\sigma_t(\theta
)}
\biggr)\Delta_{2,t}(\theta),\nonumber
\end{eqnarray}
and ${\cal R}_n$ is a reminder which is displayed below.
As in the proof of Lemma~\ref{lemderivalphabeta}, it can be
seen that
\[
S_{n,f}(\theta_0)=\frac{-1}{\delta\sqrt{n}}\sum
_{t=1}^ng_1(\eta _t)\,d_t(
\vartheta_0)+o_P(1),\qquad d_t(\vartheta)=
\pmatrix{0
\cr
d_t^{\alpha_+}
\cr
d_t^{\alpha_-}
\cr
d_t^{\beta}}.
\]
%
Using (\ref{hypo2surf}), it is easy to see that $Eg_1(\eta_1)=0$,
and thus $Eg_1^2(\eta_1)=\iota_f$. The Lindeberg central limit
theorem for martingale differences then shows that
%
\begin{equation}
\label{res2proofLAN} S_{n,f}(\theta_0) 
\stackrel{d} {\longrightarrow} {\cal N}
(0,\mathfrak{I}_f ).
\end{equation}
%
Turning to the second term of (\ref{res1proofLAN}) we first note
that, similar to (\ref{derminusdt}),
%
\[
\biggl\llvert \frac{1}{h_t}\,\frac{\partial\sigma^{\delta}_t(\theta
_0)}{\partial\theta}-d_t(
\vartheta_0)\biggr\rrvert \to0 \qquad\mbox{in $L^2$ as } t\to
\infty.
\]
Moreover, integrations by parts show that, under (\ref{hypo2surf}),
$\int y^2f''(y)\,dy=\break-2 \int y\*f'(y)\,dy=2$. It follows that
$Eg_2(\eta_1)=-\iota_{f}$. We thus have, using\break $Eg_1(\eta_1)=0$,
\begin{eqnarray*}
\mathfrak{I}_n(\theta_0)&=&
\frac{1}{n}\sum_{t=1}^n
\frac{-g_2(\eta_t)}{\delta^2} \,d_t(\vartheta_0)\,d_t'(
\vartheta_0)\\
&&{}+o_{P_{\theta_0}}(1)\to \mathfrak{I}_{f}
\qquad\mbox{in probability as } n\to\infty.
\end{eqnarray*}
Next, it can be shown that, as $n\to\infty$,
%
\begin{equation}
\label{res3proofLAN} \bigl\llVert \mathfrak{I}_n\bigl(
\theta_n^*\bigr)-\mathfrak{I}_n(\theta_0)
\bigr\rrVert \to 0 \qquad\mbox{in probability.}
\end{equation}
Finally, we show
the convergence in probability to zero of
\[
{\cal R}_n=\upsilon_n\sum_{t=1}^ng_1(
\eta_t) \frac{1}{\delta
h_t}\,\frac{\partial\sigma_t^\delta(\theta_0)}{\partial\omega} -\upsilon_n
\sqrt{n}\btau'\mathfrak{I}_n\bigl(\theta_n^*
\bigr)\be_1' -\frac12n\upsilon_n^2
\be_1\mathfrak{I}_n\bigl(\theta_n^*\bigr)
\be_1'.
\]
Noting that ${\partial
\sigma_t^\delta(\theta_0)}/{\partial\omega}$ is constant and that
$1/h_t$ converges to 0 in $L^2$ by Proposition \ref{convtoinfty},
the first term in the right-hand side
converges to zero in probability. The two other terms can be handled
similarly. The conclusion then follows from
(\ref{res1proofLAN})--(\ref{res3proofLAN}).
\end{pf*}

\begin{pf*}{Proof of Proposition \ref{LAP2tests}}
For simplicity, write $P$ instead of $P_{n,0}$. In the proof of
Theorem \ref{jointdistribution} we have seen that
\[
T_n=\frac{1}{\sqrt{n}}\sum_{t=1}^n
\frac{u_t}{\sigma_u}+o_P(1).
\]
By (\ref{lan}) and (\ref{Sn}), it follows that under $P$
\[
\pmatrix{T_n
\cr
\Lambda_{n,f}(
\theta_0+\btau/\sqrt{n},\theta_0)} \stackrel{d} {
\longrightarrow} {\cal N} \left\{\pmatrix{ 0
\vspace*{2pt}\cr
-\dfrac{\iota_{f}}{8}\tilde{
\btau}'{\cal I}\tilde {\btau}}, \pmatrix{1&c
\vspace*{2pt}\cr
c&\dfrac{\iota_{f}}{4}
\tilde{\btau}'{\cal I}\tilde{\btau}} \right\},
\]
where $\tilde{\btau}'=(\tau_2,\tau_3,\tau_4)$,
$c=-\frac{\btau'Ed_1(\vartheta_0)}{\delta\sigma_u}Eu_1g_1(\eta
_1)=c_f(\theta_0)$.
Le Cam's third lemma [see, e.g., \citet{van98}, page 90]
shows that
\[
T_n\stackrel{d} {\longrightarrow} {\cal N} \bigl(c_f(
\theta_0),1 \bigr)\qquad \mbox{under }P_{n,\btau}.
\]
The conclusion easily follows.\vadjust{\goodbreak}
\end{pf*}

\begin{pf*}{Proof of Proposition \ref{LAPtestsuralpha}} First
consider the case $\gamma_0\geq0$. In the proof of
(\ref{normalitydeleQMVdealpha}) it has been shown that
\[
\sqrt{n}(\hat{\vartheta}_n-\vartheta_0)=-
\frac{2}{\delta}{\cal I}^{-1}\frac{1}{\sqrt{n}}\sum
_{t=1}^n\bigl(1-\eta_t^2
\bigr)\,d_t+o_P(1).
\]
Moreover
\[
\Lambda_{n,f}(\theta_0+\btau/\sqrt{n},
\theta_0)=-\frac{1}{\delta
\sqrt{n}}\sum_{t=1}^n
\biggl\{1+\eta_t\frac{f'(\eta_t)}{f(\eta
_t)} \biggr\}\tilde{
\btau}'d_t- \frac{\iota_f}{8}\tilde{
\btau}'{\cal I}\tilde{\btau}+o_P(1)
\]
with $\tilde{\btau}'=(\tau_2,\tau_3,\tau_4)$. Note also that, since
$E\eta_1^4<\infty$ implies $y^3f(y)\to0$ as $|y|\to\infty$, we have
%
\begin{equation}
\label{vaut2} E\bigl(1-\eta_t^2\bigr) \biggl\{1+
\eta_t\frac{f'(\eta_t)}{f(\eta_t)} \biggr\}=2.
\end{equation}
It follows that under
$P_{n,0}^{\mathrm{S}}$
\[
\pmatrix{\sqrt{n}(\hat{\vartheta}_n-\vartheta_0)
\cr
\displaystyle \Lambda_{n,f} \biggl(\theta_0+\frac{\btau}{\sqrt{n}},\theta
_0 \biggr)} \stackrel{d} {\longrightarrow} {\cal N} \left\{
\pmatrix{0_3
\vspace*{2pt}\cr
\displaystyle \frac{-\iota_{f}}{8}\tilde{\btau}'{\cal I}
\tilde{\btau}}, \pmatrix{\displaystyle (\kappa_{\eta}-1){\cal I}^{-1}&\tilde{
\btau}
\cr
\tilde{\btau}'&\displaystyle \frac{\iota_{f}}{4} \tilde{
\btau}'{\cal I}\tilde {\btau}} \right\}.
\]
Le Cam's third lemma [see, e.g., \citet{van98}, page 90]
shows that
\[
\sqrt{n}(\hat{\vartheta}_n-\vartheta_0)\stackrel {d} {
\longrightarrow} {\cal N} \bigl(\tilde{\btau},(\kappa_{\eta}-1){\cal
I}^{-1} \bigr)\qquad \mbox{under }P_{n,\btau}^{\mathrm{S}}.
\]
We thus have shown that, in the case
$\gamma_0>0$, $\hat{\vartheta}_n$ is a regular estimator of
$\vartheta_0$, in the
sense that $\sqrt{n} (\hat{\vartheta}_n-\vartheta_0-\tilde
{\btau}/\sqrt{n} )$
converges to a distribution which does not depend on $\tilde{\btau}$. More
precisely
%
\begin{equation}
\label{regulier} \sqrt{n} (\hat{\vartheta}_n-\vartheta_0-
\tilde{\btau}/\sqrt {n} )\stackrel{d} {\longrightarrow} {\cal N} \bigl(0,(
\kappa_{\eta}-1){\cal I}^{-1} \bigr)\qquad \mbox{under
}P_{n,\btau}^{\mathrm{S}}.
\end{equation}
When $\gamma_0\leq0$, the same arguments show that $\hat{\theta}_n$
is a regular estimator of $\theta_0$
\[
\sqrt{n}(\hat{\theta}_n-\theta_0-\btau/\sqrt{n})
\stackrel {d} {\longrightarrow} {\cal N} \bigl(0,(\kappa_{\eta}-1){\cal
J}^{-1} \bigr)\qquad\mbox{under }P_{n,\btau}^{\mathrm{S}}.
\]
In the case $\gamma_0\leq0$, we thus have (\ref{regulier}) with
${\cal I}$ replaced by ${\cal I_*}$.
Now, noting that
$T_n^{\mathrm{S}}=\frac{\be'\sqrt{n}(\hat{\vartheta}_n-\vartheta
_0)}{\hat{\sigma}_{T^\mathrm{S}}}$,
and by the same arguments, it follows that
$T_n^{\mathrm{S}}\stackrel{d}{\longrightarrow}{\cal
N} (0,1 )$, under $P_{n,0}^{\mathrm{S}}$
and more generally
$T_n^{\mathrm{S}}\stackrel{d}{\longrightarrow}{\cal
N} (c_{\btau},1 )$, under $P_{n,\btau}^{\mathrm{S}}$,
where $c_{\btau}=(0,1,-1,0)\btau/\sigma_{T^\mathrm{S}}$. The
conclusion easily follows.
\end{pf*}

\begin{pf*}{Proof of Proposition \ref{coro2}} Recall that we
assume $\gamma_0\geq0$. The case $\gamma_0<0$ is obtained
similarly, replacing ${\cal I}$ by ${\cal I}_*$.
In view of Proposition \ref{LAPtestsuralpha} and
(\ref{bornepower}), the ${\mathrm C}^{\mathrm{S}}$-test is
asymptotically locally UMPU if and only if
$c_{\be'\tilde{\btau}}=\be'\tilde{\btau}/\sigma_{T^\mathrm{S}}$,
which is equivalent to $(\kappa_{\eta}-1)\iota_f=4$. By Corollary 1
in \citet{FraZak06}, the solutions of this equation are
given by (\ref{casqmleefficient}).
\end{pf*}

\section{Concluding remarks}\label{conclu}
Our framework covers the most widely used GARCH models in financial
applications. Strictly stationary models are a special case, but
symmetry tests and asymptotically valid confidence intervals for
the parameters (except the intercept) can be built without this
assumption. Surprisingly, while the asymptotic covariance matrix of
the estimators is sensitive to the stationarity of the underlying
process, an estimator which converges to the appropriate covariance
matrix in every situation can be built. Nevertheless, if the
interest is on the whole parameter vector, including the intercept,
it is important to know whether the observations come from a
stationary process or not. To this aim we derived strict
stationarity/nonstationarity tests which are very easy to
implement.

Are our results extendable to higher-order models? It seems likely
that for particular extensions involving \textit{univariate} stochastic
recurrence equations for the volatility, the asymptotic theory
derived in this paper can also be established. One key problem, to
show consistency, is to find stationary approximations to $\epsilon
_{t-j}^2/h_t$ for $j=1, 2, \ldots\,$.
For an ARCH-type model of order $q$ it suffices to take $j\leq q$.
Consider standard symmetric GARCH
models for simplicity. In the $\operatorname{GARCH}(1,1)$ case, the problem can be
circumvented because
\[
\frac{\epsilon_{t-j}^2}{h_t}=\frac{h_{t-1}}{h_t}\cdots \frac{h_{t-j}}{h_{t-j+1}}
\eta_{t-j}^2
\]
can be approximated by a
stationary process, in view of
\[
\frac{h_{t-i}}{h_{t-i+1}}\approx\frac1{\alpha \eta_{t-i}^2+\beta}
\qquad\mbox{for large $t$}.
\]

To have a glimpse of the considerable difficulties encountered when
the orders increase, consider a standard ARCH(2) model
\[
\epsilon_t=\sqrt{h_t}\eta_t,\qquad
h_t=\omega+\alpha_1\epsilon _{t-1}^2+
\alpha_2\epsilon_{t-2}^2.
\]
%
We have, neglecting $\omega$ and for $t$ large enough
${h_{t}}/{\epsilon_{t-1}^2}\approx X_t$ and
${h_{t}}/{\epsilon_{t-2}^2}\approx Y_t$ where
\[
X_t= \alpha_1 +\frac{\alpha_2}{X_{t-1}}\frac{1}{\eta_{t-1}^2},\qquad
Y_t=\alpha_2+\alpha_1\eta_{t-1}^2X_{t-1}.
\]
It is not difficult to show that the first stochastic recurrence
equation admits a strictly stationary solution $(X_t)$ under mild
assumptions on the density of $\eta_t$, whatever the values of
$\alpha_1$ and $\alpha_2$. From this solution we deduce a strictly
stationary solution $(Y_t)$ to the second equation. We thus believe
that, at least for the consistency, the ARCH(2) model is amenable to
a treatment similar to that developed in this paper, but at the
price of increasing technical difficulties. To summarize, the ratio
$h_t/h_{t-1}$ is, for large~$t$, close to (i)~a~constant in the
ARCH(1) case, (ii) an i.i.d. process in the $\operatorname{GARCH}(1,1)$ case and (iii)~the
stationary solution of a nonlinear times series model in the ARCH(2)
case. Whether or not this approach based on the resolution of
nonlinear stochastic recurrence equations could be extended is left
for further investigation.

\section*{Acknowledgments}

We are most thankful to the Editor and to three referees for their
constructive comments and suggestions.
We are also grateful to the Agence Nationale de la Recherche (ANR).

\begin{supplement}
\stitle{Supplement to ``Inference in nonstationary asymmetric GARCH models.''}
\slink[doi]{10.1214/13-AOS1132SUPP} 
\sdatatype{.pdf}
\sfilename{aos1132\_supp.pdf}
\sdescription{The
supplementary file contains an illustration concerning the optimality
of the asymmetry test, a~Monte Carlo study of finite sample
performance, an application to real time series, an explicit expression
for the matrix ${\cal I}$ in Theorem \ref{theoconsistenceJensenRahbek},
the proofs of Theorems \ref{propo1} and \ref{theoconsistenceJensenRahbekWithDelta}.}
\end{supplement}


\printaddresses

\end{document}